\documentclass[12pt]{article}

\usepackage{amsmath,amstext,amsfonts,amsthm,amssymb,eucal,diagrams}

\newcommand{\del}{\partial}

\newcommand{\ch}{\operatorname{ch}}

\renewcommand{\mod}{\operatorname{mod}}
\newcommand{\Com}{\operatorname{Com}}

\newcommand{\Tr}{\operatorname{Tr}}
\newcommand{\gr}{\operatorname{gr}}

\newcommand{\Tor}{\operatorname{Tor}}

\newcommand{\coker}{\operatorname{coker}}

\newcommand{\DD}{{\cal D}}

\newcommand{\JJ}{{\cal J}}

\newcommand{\bL}{{\bf L}}

\newcommand{\BB}{{\cal B}}

\newcommand{\mg}{{\frak m}}

\newcommand{\lan}{\langle}
\newcommand{\ran}{\rangle}

\newcommand{\CC}{{\cal C}}

\newcommand{\Spec}{\operatorname{Spec}}

\newcommand{\si}{\sigma}

\newcommand{\ga}{\gamma}
\newcommand{\de}{\delta}
\newcommand{\eps}{\epsilon}

\renewcommand{\ker}{\operatorname{ker}}
\newcommand{\im}{\operatorname{im}}

\numberwithin{equation}{section}

\newtheorem{theor}{Theorem}[subsection]
\newtheorem{thm}[theor]{Theorem}
\newtheorem{lem}[theor]{Lemma}

\newtheorem{prop}[theor]{Proposition}

\newtheorem{cor}[theor]{Corollary}
{  \theoremstyle{definition}
           \newtheorem{defn}[theor]{Definition}
           \newtheorem{rem}[theor]{Remark}
          \newtheorem{ex}[theor]{Example}
}

\newcommand{\Pf}{\noindent {\it Proof}}
\newcommand{\id}{\operatorname{id}}

\newcommand{\we}{\wedge}

\renewcommand{\AA}{{\cal A}}

\newcommand{\FF}{{\cal F}}

\newcommand{\Om}{\Omega}

\newcommand{\Hom}{\operatorname{Hom}}

\newcommand{\End}{\operatorname{End}}
\newcommand{\Res}{\operatorname{Res}}

\newcommand{\sTr}{\operatorname{str}}

\renewcommand{\a}{\alpha}
\renewcommand{\b}{\beta}
\newcommand{\om}{\omega}
\newcommand{\De}{\Delta}

\newcommand{\C}{{\Bbb C}}

\newcommand{\Z}{{\Bbb Z}}

\newcommand{\Ga}{\Gamma}

\newcommand{\wt}{\widetilde}
\newcommand{\ot}{\otimes}

\newcommand{\sub}{\subset}
\newcommand{\ed}{\qed\vspace{3mm}}

\newcommand{\tr}{\operatorname{tr}}

\newcommand{\pa}{\partial}

\newcommand{\per}{\operatorname{Per}}

\newcommand{\by}{{\bf y}}

\newcommand{\Zt}{{\Z/2}}

\newcommand{\mf}{matrix factorization}
\newcommand{\zer}{\mathbf{0}}
\newcommand{\one}{\mathbf{1}}
\newcommand{\st}{\mathrm{st}}
\renewcommand{\d}{\delta}

\newcommand{\w}{{\boldsymbol w}} 
\newcommand{\m}{\mathfrak{m}} 
\newcommand{\MF}{\mathrm{MF}} 
\newcommand{\MFt}{\mathrm{HMF}} 
\newcommand{\HMF}{\mathrm{HMF}} 
\newcommand{\Db}{\mathrm{D^b}}
\newcommand{\Dbs}{\underline{\Db}}
\newcommand{\Dbp}{\mathrm{D^b_{per}}}
\newcommand{\lrarrow}[2] {\stackrel{{\rTo^{#1}}}{{\lTo_{#2}}}}
\newcommand{\MCM}{\mathrm{MCM}}
\newcommand{\MCMs}{\underline{\MCM}}
\newcommand{\Mod}{\operatorname{Mod}}

\newcommand{\Cok}{\mathrm{coker}}
\newcommand{\Span}{\operatorname{Span}}

\newcommand{\CM}{Cohen-Macaulay}
\newcommand{\Homb}{\mathcal{H}om} 
\newcommand{\be}{\begin{equation}}
\newcommand{\ee}{\end{equation}}
\newcommand{\HRR}{Hirzebruch-Riemann-Roch}
\newcommand{\dL}{\mathbb{L}}
\newcommand{\bde}{{\bf \delta}}

\newcommand{\Hess}{\operatorname{Hess}}
\newcommand{\bx}{{\bf x}}

\title{Chern characters and \HRR\ formula for \mf s}
\date{}
\author{Alexander Polishchuk \and Arkady Vaintrob}

\begin{document}
\maketitle
\begin{abstract} 
We study the category of \mf s for an isolated hypersurface singularity.
We compute 
the canonical bilinear form on the Hochschild homology of this category.
We find explicit expressions for the Chern
character and the boundary-bulk maps 
and derive an analog of the \HRR\ formula for the Euler characteristic
of the $\Hom$-space between a pair of \mf s. We also establish $G$-equivariant
versions of these results.
\end{abstract}

\section*{Introduction}

Let $\w$ be an element of a commutative ring $R$. A {\it \mf\ of $\w$} is 
a $\Z/2$-graded finitely generated projective $R$-module $E=E_0\oplus E_1$ together with
an odd endomorphism $\d_E$ such that $\d_E^2=\w\cdot\id_E$.
Matrix factorizations have been a classical tool in the study of hypersurface singularity
algebras since the work of Eisenbud \cite{Eis}. In the geometric context 
the category of \mf s measures the failure of every coherent sheaf on the hypersurface $\w=0$
to have a finite 
locally free resolution (see \cite{Orlov}). Matrix factorizations also appear prominently
in the work of Khovanov
and Rozansky \cite{KhR} on link homology. The categories of \mf s, following the suggestion of
Kontsevich, have been used by physicists \cite{KL1,KL2} to describe 
D-branes in topological Landau-Ginzburg models. 

This paper is motivated by the desire to understand the rich structure arising
on the Hochschild homology of the category $\MF(\w)$ of \mf s of an isolated
hypersurface singularity $\w=0$, where $\w(x_1,\ldots,x_n)$ is a formal power series. 
In physics literature \cite{KR, W} the Hochschild (co)homology of this category has an interpretation
as the state space for the closed string sector of the open-closed topological string theory
associated with the Landau-Ginzburg model of the potential $\w$.
The rigorous computation of this Hochschild homology has been done recently by Dyckerhoff
\cite{Dyck}. In our paper we derive explicit
formulas for some of the natural structures on this space using the tools of the theory of dg-categories 
\cite{Keller-dg, Toen, TV}.
The Hochschild homology of the category of \mf s of the quasihomogeneous isolated singularity
(in the orbifold setting) can be identified with the state space of the cohomological field theory (in
the sense of \cite{KM}) constructed by Fan, Jarvis and Ruan in \cite{FJR}.
The results of the present paper
are used in the sequel \cite{PV} to construct a purely algebraic version of the Fan-Jarvis-Ruan
theory.

The category $\MF(\w)$ for an isolated singularity 
fits naturally into the framework of noncommutative geometry developed from the
point of view of dg-categories or $A_{\infty}$-categories (see \cite{KKP}, \cite{Dyck}). 
As shown in \cite{Dyck} it provides an example of
a smooth and proper noncommutative space (in the sense of \cite{KS}).
The classical \HRR\ formula for coherent sheaves on smooth projective varieties
was recently generalized by Shklyarov \cite{Shk} to such noncommutative spaces
(see also \cite{CalW} where similar ideas are developed in the classical case).
He showed that the Hochschild homology $HH_*(\CC)$
of a smooth proper dg-category $\CC$ is
equipped with a canonical nondegenerate bilinear form $\lan\cdot,\cdot\ran$. 
His categorical \HRR\ formula 
expresses the Euler characteristic of the $\Hom$-spaces between two objects in the 
derived category of $\CC$ in terms of Chern characters (aka Euler classes) taking values in $HH_*(\CC)$
and the form $\lan\cdot,\cdot\ran$. 
In fact, there exists an even more general formula
computing the traces of certain endomorphisms of the $\Hom$-spaces between two objects
(see \cite[Thm. 16]{CalW} and Theorem \ref{cardy-prop} below).
In the case of a Calabi-Yau category $\CC$ this generalized formula is equivalent to the
Cardy condition for the corresponding open-closed 2d TQFT (see \cite[Thm. 15]{CalW}).

In concrete situations the difficulty shifts to
calculating explicitly the Hochschild homology of the corresponding category
along with the Chern character map and the canonical bilinear form.
In our paper we work out these ingredients in the case of the 
($\Zt$-graded) dg-category of \mf s of an isolated hypersurface singularity $\w$, including
the $G$-equivariant version, where $G$ is a finite group of symmetries of $\w$. This leads to a concrete \HRR\ formula for \mf s.

Now we will give an explicit formulation of our results (in the non-equivariant setting).
By the result of Dyckerhoff \cite{Dyck}, for an isolated singularity $\w\in k[[x_1,\ldots,x_n]]$
the Hochschild homology of the ($\Zt$-graded dg-)category $\MF(\w)$ of \mf s of $\w$ is
isomorphic as a $\Zt$-graded vector space to the Milnor ring of $\w$ (up to a shift of grading):
\be\label{hoch-milnor-eq}
HH_*(\MF(\w))\simeq \AA_{\w}\cdot d\bx[n],
\ee
where $\AA_{\w}=k[[x_1,\ldots,x_n]]/\JJ_{\w}$ with
$\JJ_{\w}=(\pa_1\w,\ldots,\pa_n\w)$ and $d\bx=dx_1\we\ldots\we dx_n$.
We will derive 
the following formula for the Chern character $\ch(\bar{E})\in HH_0(\MF(\w))$ 
of a \mf\ $\bar{E}=(E,\d_E)$ (see Theorem \ref{ch-thm}):
\begin{equation}\label{ch-eq}
\ch(\bar{E})=\sTr_R(\pa_n \d_E\cdot\ldots\cdot\pa_1 \d_E)\cdot d\bx \ \mod\ 
\JJ_{\w}\cdot d\bx,
\end{equation}
where we view $\d_E$ as a matrix
with entries in $R=k[[x_1,\ldots,x_n]]$ (by choosing a basis of the free module $E$)
and take partial derivatives $\pa_i=\pa/\pa x_i$ component-wise, and $\sTr_R$ is the supertrace
of a matrix with entries in $R$.
More generally, 
there is a canonical map 
\begin{equation}\label{bb-map-eq}
\tau^{\bar{E}}:\Hom^*(\bar{E},\bar{E})\to HH_*(\MF(\w))
\end{equation}
(called the ``boundary-bulk map" in the context of topological strings),
such that $\ch(\bar{E})=\tau^{\bar{E}}(\id_E)$ and the formula \eqref{ch-eq} generalizes to
\begin{equation}\label{iota-eq}
\tau^{\bar{E}}(\a)=\sTr_R(\pa_n \d_E\cdot\ldots\cdot\pa_1 \d_E\circ\a)\cdot d\bx \ \mod\ 
\JJ_{\w}\cdot d\bx
\end{equation}
where $\a$ is an arbitrary endomorphism of $\bar{E}$.

We will also show in Corollary \ref{form-cor} that the formula \eqref{ch-eq} 
leads to the identification of the canonical bilinear form on $HH_*(\MF(\w))$
with the form
\begin{equation}\label{form-eq}
\lan f\ot d\bx, g\ot d\bx\ran=(-1)^{{n\choose 2}}\Res(f\cdot g),
\end{equation}
where $\Res$ is 
the linear functional on the Milnor ring $\AA_\w$ given by the 
generalized residue:
$$\Res(f)=\Res_{k[x]/k}\begin{bmatrix} f(x)\cdot dx_1\we\ldots\we dx_n \\ \pa_1\w,\ldots,\pa_n\w\end{bmatrix}
$$
(see \cite[III.9]{Hart-RD}, \cite[ch.\ V]{GH}).
As a consequence we obtain an analog of the \HRR\ formula for the Euler characteristic of the
$\Zt$-graded space $\Hom^*(\bar{E},\bar{F})$: 
\begin{equation}\label{HRR-eq}
\chi(\bar{E},\bar{F})=\lan\ch(\bar{E}),\ch(\bar{F})\ran.
\end{equation}
More generally, for $\a\in\Hom^*(\bar{E},\bar{E})$ and $\b\in\Hom^*(\bar{F},\bar{F})$ 
we have
\begin{equation}\label{cardy-eq-intro}
\sTr_k(m_{\a,\b})=\lan\tau^{\bar{E}}(\a),\tau^{\bar{F}}(\b)\ran,
\end{equation}
where $m_{\a,\b}$ is the endomorphism of $\Hom^*(\bar{E},\bar{F})$ sending
$f$ to $(-1)^{|\a|\cdot|\b|+|\a|\cdot|f|}\b\circ f\circ\a$
(see Theorem \ref{HRR-thm}).
All terms in the right-hand sides of \eqref{HRR-eq} and \eqref{cardy-eq-intro} can be
explicitly expressed via \eqref{iota-eq} and \eqref{form-eq}
in terms of partial derivatives of $\d_E$, $\d_F$ and $\w$. 

In the case when $n$ is odd the formula \eqref{HRR-eq} implies that $\chi(\bar{E},\bar{F})$
is identically zero, since in this case $HH_0(\MF(\w))=0$. This was conjectured by Hailong Dao
(see \cite[Conj.\ 3.15]{Dao} and Remark \ref{odd-rem} below).

Note that in some particular cases the \HRR\ formula for \mf s was proved in \cite{W}.
Our formula \eqref{iota-eq}
for the map $\tau^{\bar{E}}$ is almost identical to the formula for the {\it boundary-bulk map}
in the corresponding Landau-Ginzburg model for open topological strings (see \cite{KL2}, \cite{Segal}):
we get an extra sign (see Corollary \ref{perm-cor}).
Similar expression also appears in the explicit version of the Serre duality for \mf s worked out
by Murfet in \cite{Murfet}. In the present paper Serre duality does not appear; 
this connection will be discussed elsewhere.

We will also establish $G$-equivariant versions of formulas \eqref{ch-eq}, \eqref{iota-eq}, 
\eqref{form-eq}, \eqref{HRR-eq} and \eqref{cardy-eq-intro}
 (see Theorems \ref{G-hoch-thm} and \ref{G-HRR-thm}).

In the case of a quasihomogeneous singularity one can also consider 
$\Z$-graded versions of the categories of \mf s (see e.g. \cite{W}, \cite{Orlov-gr}).
An analog of the \HRR\ formula for these categories follows 
from this
formula for the category of $G$-equivariant \mf s, where $G$ is an 
appropriate cyclic group (see section \ref{graded-sec}).

The paper is organized as follows. In section \ref{hoch-gen-sec}  we review some general
facts about Hochschild homology of dg-categories. In section \ref{functoriality-sec}
we give a slightly nonstandard construction of the maps on Hochschild homology induced
by dg-functors. This construction fits well with our method of computing the boundary-bulk maps
(and can be shown to be equivalent to the standard one). 
We also present in \ref{cardy-sec}
a generalized categorical version of the \HRR\ theorem (a version of
the Cardy condition discussed in \cite{CalW}). In section \ref{sec:mf} we collect useful
(and mostly well known) results on \mf s, including calculation of the Hochschild homology of
the dg-category $MF(\w)$ by the techniques
of Dyckerhoff~\cite{Dyck} and present a $G$-equivariant version of this calculation. 
Section \ref{Chern-sec} contains the main computation leading to the explicit formula for the Chern
character. Finally, in section \ref{form-sec} this formula is used to find the expression for the canonical bilinear form on Hochschild homology (using calculations from \cite{KR}) and to derive a \HRR\ formula
for matrix factorizations.
We also present in \ref{k-st-sec} the explicit calculation of the boundary-bulk map \eqref{bb-map-eq}
for the special matrix factorization $k^{\st}$ (the {\it stabilization of the residue field}).

\bigskip

\noindent
{\it Notations and conventions.} We work over a fixed ground field $k$ of characteristic zero. 
All the dg-categories considered in this paper are assumed to be $k$-linear. We denote
the tensor product of $k$-vector spaces simply by $V\otimes W$. For a $\Zt$-graded vector space
$V=V^\zer\oplus V^\one$ we denote by $J_V$ the {\it grading operator} that sends a
homogeneous vector $v$ to $(-1)^{|v|}\cdot v$.

\vspace{2mm}

\noindent
{\it Acknowledgments}. We thank the referees for many useful comments and suggestions. We are especially
grateful to one of the referees who provided a simpler proof of Proposition \ref{res-prop}.
Both authors would like to thank the IHES, where this work was done,
for hospitality and stimulating atmosphere. The first author was partially supported by the NSF
grant DMS-0601034.

\section{Hochschild homology for smooth proper dg-categories}
\label{hoch-gen-sec}

Here we review the definition and some properties of the Hochschild homology for differential graded categories (dg-categories).
Nice surveys of the relevant facts can be found in the
papers \cite{Keller-dg} and \cite{Toen-dg}. An extension of these techniques to the $\Z/2$-graded case 
is explained in \cite{Dyck}. We always assume our dg-categories to be small (or quasi-equivalent
to such).

  
\subsection{Modules and Hochschild homology for dg-categories}

The notion of a $k$-linear dg-category is a generalization of the notion of dg-algebra (which
is a dg-category with one object), 
and most constructions for dg-algebras can be similarly defined for dg-categories.
For instance, a (left) module of a dg-category $\CC$ is a dg-functor $\CC\to\Com(k)$,
where $\Com(k)$ is the category of complexes over $k$. We denote by $\CC-\mod$ the dg-category
of $\CC$-modules.
For every dg-category $\CC$ we denote by $\CC^{op}$ the opposite dg-category with the
same set of objects but with the composition $f\circ g$ replaced with $(-1)^{|f|\cdot|g|}g\circ f$.
For an object $E\in\CC$ we denote the same object viewed as an object of $\CC^{op}$
by $E^{\vee}$. We define right $\CC$-modules as modules over $\CC^{op}$.
For a pair of dg-categories $\CC$ and $\CC'$ we define $\CC-\CC'$-bimodules as
right modules over $\CC^{op}\ot \CC'$. Similarly to the ordinary Morita theory we can use
the tensor product with a $\CC-\CC'$-bimodule $X$ 
to obtain a dg-functor $T_X$ from the category of right $\CC$-modules to that
of right $\CC'$-modules (see \cite[Sec.\ 6.1]{Keller-derived}). 
The derived category $D(\CC)$ of a dg-category $\CC$ is defined as the localization of
the category of right $\CC$-modules with respect to quasi-isomorphisms (see \cite[Sec.\ 3]{Keller-dg}).
The above tensor product functor $T_X$ gives rise to the left derived functor 
\begin{equation}\label{tensor-kernel-eq}
\bL T_X:D(\CC)\to D(\CC'):M\mapsto M\ot^{\dL}_{\CC} X.
\end{equation}
By analogy with the theory of Fourier-Mukai transforms
we say that $\bL T_X$ is the functor associated with the kernel $X$.
In particular, the identity functor $D(\CC)\to D(\CC)$ is represented by the
"diagonal" $\CC-\CC$-bimodule $\De=\De_{\CC}$ given by
the dg-functor 
$$(\CC^{op}\ot\CC)^{op}=\CC\ot\CC^{op}\to \Com(k): E\ot F^{\vee}\mapsto\Hom_{\CC}(F,E).$$

The Hochschild homology $HH_*(\CC)$ of a (small) dg-category $\CC$ is usually defined
in terms of the Hochschild chain complex (see e.g., \cite[Sec.\ 5.3]{Keller-dg}).
As in the case of dg-algebras, this is equivalent to defining $HH_*(\CC)$ as the derived tensor
product of the $\CC-\CC$-bimodule $\De$ with itself. Indeed, the Hochschild chain complex 
arises when one computes the derived tensor product 
$\De \ot^{\dL}_{\CC^{op}\ot\CC} \De$ using the bar resolution of $\De$ 
in the category of $\CC-\CC$-bimodules (see e.g. \cite[Sec.\ 6.6]{Keller-derived}).
Following Toen (see \cite[Sec.\ 6.3]{Dyck} and \cite[Sec.\ 5.2.3]{Toen-dg}) let us introduce
the functor
\be\label{Tr-def-eq}
\Tr_{\CC}: D(\CC^{op}\ot\CC)\to D(k):M\mapsto M \ot^{\dL}_{\CC^{op}\ot\CC}\De
\ee
of the derived tensor product with $\De$.
Then the definition of the Hochschild homology can be 
rewritten as follows:
\begin{equation}\label{hochschild-hom-def}
HH_*(\CC)=\De\ot^{\dL}_{\CC^{op}\ot\CC} \De=\Tr_{\CC}(\De).
\end{equation}

\begin{rem} The notation $\Tr_{\CC}$ comes from the analogy with the computation of the trace of an
integral transform through the restriction of its kernel to the diagonal. For the simplest
example showing the ``trace-like" nature of this functor, take $\CC$ to be the algebra
$k^I$, the direct sum of the ground field $k$ over a finite index set $I$. Then
$\CC-\CC$-bimodules correspond to $I\times I$-graded vector spaces $(V_{ij})_{i,j\in I}$, and the functor
$\Tr_{\CC}$ sends $(V_{ij})$ to the vector space $\bigoplus_{i\in I}V_{ii}$. Note however that
the functor $\Tr_{\CC}$ is different from the categorical trace introduced in the context of
$2$-categories by Ganter
and Kapranov in \cite{GK}. Instead of the functor of tensor product with $\De$ 
they use the functor $\Hom_{\CC}(\De,?)$ 
(which leads to the definition of the Hochschild cohomology).
\end{rem}



For any dg-category $\AA$ we denote by $H^0\AA$ the corresponding {\it homotopy category}
obtained by passing to $0$th cohomology of the $\Hom$-spaces.
Recall that we have the Yoneda embedding $H^0\AA\to D(\AA)$
sending $A\in\AA$ to the representable dg-functor 
$$h_A=\Hom_{\AA}(?,A):\AA^{op}\to\Com(k).$$
Note that for any left $\AA$-module $M$ and any $A\in\AA$
one has
\begin{equation}\label{M(A)-eq}
h_A\ot^{\dL}_{\AA} M= h_A\ot_{\AA} M= M(A).
\end{equation}
In particular, for $F,E\in\CC$ we have
$$\Tr_{\CC}(h_{F^\vee\ot E})=h_{F^\vee\ot E}\ot_{\CC^{op}\ot\CC}\De=\Hom_{\CC}(F,E).$$

\begin{rem}\label{Tr-sign-rem} 
When computing the action of $\Tr_{\CC}$ on morphisms some signs appear due to
the standard sign convention. Namely, for a pair of morphisms $e:E_1\to E_2$,
$f:F_2\to F_1$ the induced morphism
$$\Tr_{\CC}(h_{F_1^{\vee}\ot E_1})=\Hom_{\CC}(F_1,E_1)\to\Hom_{\CC}(F_2,E_2)=
\Tr_{\CC}(h_{F_2^{\vee}\ot E_2})$$
is given by 
$$x\mapsto (f^{\vee}\ot e)(x)=(-1)^{|e|\cdot |f|+|x|\cdot |f|}\cdot e\circ x\circ f.$$
\end{rem}

Note that under the natural equivalence $\si:D(\CC\ot\CC^{op})\simeq D(\CC^{op}\ot\CC)$ 
we have isomorphisms 
$$\si(\De_{\CC^{op}})\simeq\De_{\CC},$$
\be\label{Tr-sigma-eq}
\Tr_{\CC}\circ\si\simeq\Tr_{\CC^{op}},
\ee
which induce an isomorphism
$$\phantom{x}^{\vee}:HH_*(\CC^{op})\wt{\to} HH_*(\CC)$$

The {\it perfect derived category}
$\per(\CC)\sub D(\CC)$ is the full subcategory of $D(\CC)$ defined as the minimal thick triangulated
subcategory containing all representable functors. It coincides with the subcategory of compact 
objects in $D(\CC)$ (see \cite[Cor.\ 3.7]{Keller-dg}).
Note that there is a natural equivalence
\begin{equation}\label{per-op-eq}
\per(\CC)^{op}\rTo{\sim}\per(\CC^{op}):E\mapsto E^\vee
\end{equation}
sending $h_C$ to $h_{C^\vee}$
that corresponds to the standard duality of left and right perfect modules (cf. \cite[(3.6)]{Shk} and
the proof of \cite[Thm.\ 4.6]{Keller-der-inv}). 

We will also need the following dg-versions 
of the derived category and the perfect subcategory. We define
$D_{dg}(\CC)$ as the dg-category of cofibrant right $\CC$-modules with respect to the natural model structure on
the category of right $\CC$-modules and closed morphisms
(in \cite{Toen} this dg-category is denoted by $\widehat{\CC}$, we use
the notation of \cite{Keller-dg}).
By \cite[Prop.\ 3.5]{Toen}, we have $D(\CC)\simeq H^0D_{dg}(\CC)$. 
We define $\per_{dg}(\CC)$ as the full dg-subcategory of $D_{dg}(\CC)$ consisting of
homotopically finitely presented $\CC^{op}$-modules. More explicitly, $\per_{dg}(\CC)$ is the closure
of the essential image of the Yoneda embedding under taking retracts, shifts and homotopy push-outs.
We have
$H^0\per_{dg}(\CC)=\per(\CC)\sub D(\CC)$ (see \cite[Sec.\ 7]{Toen}).
For $E,F\in\per_{dg}(\CC)$ we denote by $\Hom_\CC(E,F)$ the morphism space
in $\per_{dg}(\CC)$.

We will need the following property of the functors $\Tr_{\CC}$.

\begin{lem}\label{Tr-lem} 
(i) For dg-categories $\CC,\DD$ and objects $F\in\per(\CC^{op}\ot \DD)$ and
$G\in\per(\DD^{op}\ot \CC)$ 
there is a canonical functorial isomorphism
\be\label{Tr-FG-eq}
\Tr_{\DD}(G\ot_{\CC} F)\simeq\Tr_{\CC}(F\ot_{\DD} G)
\ee
in $D(k)$.

\noindent
(ii) For $E\in\per_{dg}(\DD)$ one has a canonical isomorphism
$$\Tr_{\DD}(E^\vee\ot E)\simeq\Hom_{\DD}(E,E),$$
 where $E^\vee(D)=\Hom_{\DD^{op}-\mod}(E,h_D)$ corresponds to $E$
by duality \eqref{per-op-eq}.
\end{lem}

\Pf . (i) Consider the tensor product $F\ot G\in\per(\CC^{op}\ot\DD\ot\DD^{op}\ot\CC)$.
We have a natural functor
$$D(\CC^{op}\ot\DD\ot\DD^{op}\ot\CC)\to D(k)$$
induced by $\Tr_{\CC}$ and $\Tr_{\DD}$.
The isomorphism \eqref{Tr-FG-eq} reflects two ways of evaluating this functor on $F\ot G$:
by either first applying $\Tr_{\CC}$ and then $\Tr_{\DD}$, or vice versa.
The key observation here is that there is a canonical isomorphism
\begin{equation}\label{F-G-tensor-diag-eq}
F \ot^{\dL}_{\DD} G\to (F\ot G)\ot^{\dL}_{\DD\ot\DD^{op}}\De_{\DD}
\end{equation}
and similarly for $G\ot_{\CC} F$. Indeed, we have a canonical morphism between the underived tensor products
for arbitrary $F$ and $G$
(not necessarily perfect), constructed using the definitions of tensor products on both sides as cokernels and
the morphism
$$\bigoplus_{D\in\DD} F(C_1,D^\vee)\ot G(D,C_2)\to\bigoplus_{D_1,D_2\in\DD}
F(C_1,D_1^\vee)\ot G(D_2,C_2^\vee)\ot \Hom_{\DD}(D_2,D_1)$$
given by the identity maps in $\Hom_{\DD}(D_2,D_1)$ for $D_2=D_1=D$. Applying this to projective resolutions
of $F$ and $G$ (see \cite[Sec.\ 3.1]{Keller-derived}) we obtain the morphism \eqref{F-G-tensor-diag-eq}.
It is enough to check that this map is an isomorphism when $F$ and $G$ are
representable. If $F=h_{C_1^{\vee}\ot D_1}$ and $G=h_{D_2^\vee\ot C_2}$ then both
parts of \eqref{F-G-tensor-diag-eq} are isomorphic to
$$\Hom_{\DD}(D_2,D_1)\ot h_{C_1^{\vee}\ot C_2}$$
as follows from \eqref{M(A)-eq}.

\noindent
(ii) Taking $\CC=k$, $F=E$ and $G=E^\vee$ in (i) we obtain an isomorphism
$$\Tr_{\DD}(E^\vee\ot E)\simeq E\ot_{\DD} E^{\vee}.$$
It remains to use the natural quasi-isomorphism 
$$M\ot_{\DD} E^\vee\rTo{\sim}\Hom_{\DD^{op}-\mod}(E,M)$$
for a $\DD^{op}$-module $M$ (see \cite[Lem.\ 6.2]{Keller-derived}).
\ed





Another basic property of Hochschild homology is the K\"unneth isomorphism (recall that we work
over a field).

\begin{prop} For dg-categories $\CC$ and $\DD$ we have a canonical isomorphism
\be\label{Kunneth}
HH_*(\CC\ot\DD)\simeq HH_*(\CC)\ot HH_*(\DD).
\ee
\end{prop}

\Pf . The standard proof for algebras using the Hochschild chain complexes 
(see \cite[Prop.\ 9.4.1]{Weibel})
can be generalized to dg-categories (see \cite[Sec.\ 2.4]{Shk}). Here we give another proof based on the definition \eqref{hochschild-hom-def}. Note that for any dg-categories $\CC_1,\CC_2$ 
we have a natural functor of external tensor product
$$D(\CC_1)\times D(\CC_2)\to D(\CC_1\ot\CC_2):(M,N)\mapsto M\ot N,$$
where
$$(M\ot N)(C_1^\vee,C_2^\vee)=M(C_1^\vee)\ot N(C_2^\vee).$$
We have 
\be\label{kunneth-diagonal}
\De_{\CC\ot\DD}\simeq \De_{\CC}\ot \De_{\DD}\in D(\CC^{op}\ot \CC\ot \DD^{op}\ot\DD)\simeq
D((\CC\ot\DD)^{op}\ot (\CC\ot\DD)).
\ee
Furthermore, for $M_1\in D(\CC)$, $M_2\in D(\CC^{op})$, 
$N_1\in D(\DD)$ and $N_2\in D(\DD^{op})$ one has a canonical isomorphism
$$(M_1\ot N_1)\ot^{\dL}_{\CC\ot\DD}(M_2\ot N_2)\simeq 
(M_1\ot^{\dL}_{\CC} M_2)\ot (N_1\ot^{\dL}_{\DD} N_2).$$
Combining this formula with \eqref{kunneth-diagonal} we obtain the required isomorphism
$$\De_{\CC\ot\DD}\ot^{\dL}_{(\CC\ot\DD)^{op}\ot(\CC\ot\DD)}\De_{\CC\ot\DD}\simeq
(\De_{\CC}\ot^{\dL}_{\CC^{op}\ot\CC}\De_{\CC})\ot (\De_{\DD}\ot^{\dL}_{\DD^{op}\ot\DD}\De_{\DD}).$$
\ed

\subsection{Functoriality}\label{functoriality-sec}


Two dg categories $\CC$ and $\DD$ are called
{\it dg Morita equivalent} if there is an equivalence of their derived categories
$D(\CC)\to D(\DD)$ given by a composition of equivalences associated with bimodules
and their inverses (see \cite[3.8]{Keller-dg}).
 

From now on we will consider only dg-categories $\CC$ that are
dg Morita equivalent to
homologically smooth and proper dg-algebras (see \cite[4.7]{Keller-dg}). In other words,
the $\CC-\CC$-bimodule $\De_{\CC}$ should be perfect, the complexes $\Hom_{\CC}(A,B)$ for $A,B\in\CC$ should have finite dimensional cohomology, and the category $D(\CC)$ should have a compact
generator.
Such dg-categories can be characterized by the condition that $\per_{dg}(\CC)$ is 
{\it saturated} (see \cite[Sec.\ 2.2]{TV}). 

An important property of these categories is that every dg-functor
$\per_{dg}(\CC)\to\per_{dg}(\CC')$ comes from 
a perfect $\CC-\CC'$-bimodule, defined uniquely up to a canonical isomorphism in the
derived category of bimodules (see \cite[Sec.\ 2.2]{TV} and
\cite[Sec.\ 5.4]{Toen-dg}). 
Specifically, a perfect $\CC-\CC'$-bimodule $X$ gives rise to the dg-functor
$T_X: M\mapsto M\ot_{\CC} X$ from $\CC^{op}$-modules to $(\CC')^{op}$-modules
that sends $\per_{dg}(\CC)$ to $\per_{dg}(\CC')$.
The induced functor $\per(\CC)\to\per(\CC')$ is the restriction of \eqref{tensor-kernel-eq}.


An example of a saturated dg-category is provided
by the bounded dg-derived category of coherent sheaves on a smooth projective variety.
As we will see below in section \ref{sec:generator}, the perfect derived category $\per_{dg}(\CC)$
of the category $\CC$ of \mf s of an isolated singularity is also saturated.



As was shown by Keller (see \cite{Keller-exact}), Hochschild homology is invariant under dg Morita equivalence. Furthermore, 
for a pair of dg-categories $\CC$, $\DD$ as above and a dg-functor
$F:\per_{dg}(\CC)\to\per_{dg}(\DD)$ 
one can define a natural map
\be\label{F*-eq}
F_*:HH_*(\CC)\to HH_*(\DD)
\ee
as follows.

Let $X$ be a perfect $\CC-\DD$-bimodule such that $F$ is given by
$T_X$, the functor of the tensor product with $X$.
Consider the $\DD-\CC$-bimodule $X^T$ given by 
$$X^T(D,C^\vee)=\Hom_{\DD^{op}-\mod}(X(C,?),h_D).$$
Let 
\begin{equation}\label{G-fun-eq}
G=T_Y:\per_{dg}(\DD)\to\per_{dg}(\CC)
\end{equation}
be the dg-functor induced by the cofibrant replacement $Y\to X^T$ in the category of
$\CC\ot\DD^{op}$-modules.
The functors induced by $F$ and $G$ on the perfect derived categories form an adjoint pair
(see \cite[Sec.\ 6.2]{Keller-derived}). We will refer to $G$ as {\it right quasi-adjoint} functor to $F$.
We are going to construct canonical {\it quasi-adjunction maps}
\begin{equation}\label{quasi-adj-maps}
\De_{\CC}\to G\circ F, \ \ F\circ G\to \De_{\DD}
\end{equation}
in $\per(\CC^{op}\ot\CC)$ and $\per(\DD^{op}\ot\DD)$, respectively, inducing the adjunction maps
for the corresponding functors on perfect derived categories.
Note that the compositions $G\circ F$ and $F\circ G$ correspond to the tensor products 
$X\ot_{\DD}Y\in \per(\CC^{op}\ot\CC)$ and
$Y\ot_{\CC}X\in \per(\DD^{op}\ot\DD)$, respectively.
Also, since the morphisms \eqref{quasi-adj-maps} we want  to define live
in the derived categories, we can replace $Y$ by $X^T$.
By definition, the $\DD-\DD$-bimodule $X^T\ot_{\CC}X$ is given by
$$(D_1,D_2^\vee)\mapsto \coker(\phi),$$
where $\phi$ is the natural map
$$\phi:\bigoplus_{C_1,C_2} X^T(D_1,C_2^\vee)\ot\Hom_{\CC}(C_1,C_2)\ot X(C_1,D_2^\vee)
\to \bigoplus_C X^T(D_1,C^\vee)\ot X(C,D_2^\vee).$$
The required morphism $Y\ot_{\CC}X\to\De_{\DD}$ is induced by the
natural evaluation map
\begin{equation}\label{XT-X-eval-map}
\begin{diagram}
X^T(D_1,C^\vee)\ot X(C,D_2^\vee)=&\Hom_{\DD^{op}-\mod}(X(C,?),h_{D_1})\ot X(C,D_2^\vee)\\
&\dTo{}\\
&h_{D_1}(D_2^\vee)&\rTo^{=}\Hom_{\DD}(D_2,D_1).
\end{diagram}
\end{equation}
To define the morphism $\De_{\CC}\to X\ot_{\DD} Y$ we use the natural quasi-isomorphism of
bimodules
$$X\ot_{\DD} X^T\rTo{\sim} Z:=\bigl((C_1,C_2^\vee)\mapsto\Hom_{\DD^{op}-\mod}(X(C_2,?),X(C_1,?))\bigr)$$
(see \cite[Lem.\ 6.2]{Keller-derived}) and observe that there is a natural map $\De_{\CC}\to Z$
given by the functoriality of $X(?,?)$ in the first argument.

Now we define the desired map \eqref{F*-eq} as the composition
$$\Tr_{\CC}(\De_{\CC})\to\Tr_{\CC}(G\circ F)\simeq\Tr_{\DD}(F\circ G)\to\Tr_{\DD}(\De_{\DD})$$
of maps induced by the morphisms \eqref{quasi-adj-maps}
and by the isomorphism of Lemma \ref{Tr-lem}(i).

\begin{lem}\label{fun-lem} 
For functors $F:\per_{dg}(\CC)\to\per_{dg}(\DD)$ and $F':\per_{dg}(\BB)\to\per_{dg}(\CC)$ 
we have
$(F\circ F')_*=F_*\circ F'_*$.
\end{lem}

\Pf . Let $G:\per_{dg}(\DD)\to\per_{dg}(\CC)$ and $G':\per_{dg}(\CC)\to\per_{dg}(\BB)$ be the corresponding right quasi-adjoint functors.
The required equality follows from the commutative diagram
\begin{diagram}
&&\Tr_{\BB}(G'F')&&&&\Tr_{\CC}(\De)&&&&\Tr_{\DD}(FG)\\
&\ruTo{}&&\rdTo{}&&\ruTo{}&&\rdTo{}&&\ruTo{}&&\rdTo{}\\
\Tr_{\BB}(\De)&&&&\Tr_{\CC}(F'G')&&&&\Tr_{\CC}(GF)&&&&\Tr_{\DD}(\De)\\
&\rdTo{}&&\ldTo{}&&\rdTo{}&&\ruTo{}&&\luTo{}&&\ruTo{}\\
&&\Tr_{\BB}(G'GFF')&&\wt{\to}&&\Tr_{\CC}(F'G'GF)&&\wt{\to}&&\Tr_{\DD}(FF'G'G)
\end{diagram}
since the composition of the upper sequence of $6$ diagonal arrows is equal to $F_*\circ F'_*$,
while the composition of the lower $4$ arrows is $(F\circ F')_*$.
\ed

\begin{rem} There are two other ways to construct the maps on Hochschild homology induced by dg-functors: 
one using explicit complexes (see e.g. \cite[Sec.\ 2.3]{Shk}) and another using Serre functors (as in \cite{CalW}).
It is possible to connect our construction to both---this will be discussed elsewhere.
\end{rem}

\begin{defn} (see \cite[Sec.\ 1.2]{Shk})
The {\it Chern character} of an object $E\in\per(\CC)$ is the element
\be\label{ch-def-eq}
\ch(E)=(\one_E)_*(1)\in HH_0(\CC),
\ee
where $\one_E:\per_{dg}(k)\to\per_{dg}(\CC)$ is the functor sending $k$ to $E$.
\end{defn}

The Chern character with values in Hochschild homology is also called the
Euler character or the Euler class (see e.g., \cite{BNT}, \cite{Keller-cyclic-spaces}).

From our assumptions that $\De_{\CC}$ is a perfect bimodule and that $\Hom$-complexes for $\CC$ are perfect it follows that the functor \eqref{Tr-def-eq} restricts to a functor
\begin{equation}\label{tr-per-eq}
\Tr_{\CC}:\per(\CC^{op}\ot \CC)\to\per(k).
\end{equation}
Similarly, we have the dg-functor
$$\Tr_{\CC}^{dg}:\per_{dg}(\CC^{op}\ot \CC)\to\per_{dg}(k).$$
Combining the functoriality for this functor with 
the K\"unneth isomorphism \eqref{Kunneth} we obtain a canonical pairing
\begin{equation}\label{pair-eq}
\lan\cdot,\cdot\ran_{\CC}:HH_*(\CC^{op})\ot HH_*(\CC)\to k
\end{equation}
(cf. \cite[Sec.\ 1.2]{Shk}).
This leads to an analog of the \HRR\ formula for the Euler characteristics of
the $\Hom$-spaces:
\begin{equation}\label{HRR-cat-eq}
\chi(\Hom_{\CC}(E,F))=\lan\ch(E^{\vee}),\ch(F)\ran_{\CC}
\end{equation}
(see \cite[(1.2)]{Shk} and Theorem \ref{cardy-prop} where a more general result is proved).
Note that from \eqref{Tr-sigma-eq} we obtain
$$\lan h,h'\ran_{\CC^{op}}=\lan h',h\ran_{\CC},$$
where $h\in HH_*(\CC)$ and $h'\in HH_*(\CC^{op})$.

The pairing \eqref{pair-eq} is nondegenerate under our assumptions on $\CC$ (see
\cite[Thm.\ 6.2]{Shk}). The proof of this fact
uses
the following connection of the canonical pairing $\lan\cdot,\cdot\ran_{\CC^{op}}$ 
with the diagonal object
$\De_{\CC}\in \per(\CC^{op}\ot\CC)$. Consider the  Chern character
$\ch(\De_{\CC})\in HH_*(\CC^{op})\ot HH_*(\CC)$.
Then the argument of \cite[Thm.\ 6.2]{Shk} (based on Lemma \ref{fun-lem} and the fact that
$\De_{\CC}$ represents the identity functor) shows that
\begin{equation}\label{inv-De-eq}
(\lan\cdot,\cdot\ran_{\CC^{op}}\ot\id)(h\ot \ch(\De_{\CC}))=h
\end{equation}
for all $h\in HH_*(\CC)$.
Similarly,
\begin{equation}\label{inv-De-eq2}
(\id\ot \lan\cdot,\cdot\ran_{\CC^{op}})(\ch(\De_{\CC})\ot h')=h'
\end{equation}
for all $h'\in HH_*(\CC^{op})$.
The equations \eqref{inv-De-eq} and \eqref{inv-De-eq2} imply that the 
bilinear form $\lan\cdot,\cdot\ran_{\CC^{op}}$ is nondegenerate and is equal to the inverse of
the tensor $\ch(\De_{\CC})$. In other words, $\ch(\De_{\CC})$ is the {\it Casimir element} corresponding
to the nondegenerate form $\lan\cdot,\cdot\ran_{\CC^{op}}$.

We will use the following way of computing the Chern character.
Note that for every $E\in\per_{dg}(\CC)$ there is a natural map 
\begin{equation}\label{cE-map}
c_E:E^{\vee}\ot E\to \De_{\CC}
\end{equation}
of $\CC-\CC$-bimodules. Indeed, this is just a particular case of the map \eqref{XT-X-eval-map}  
(with $\CC=k$ and $\DD=\CC$).

\begin{prop}\label{ch-lem} 
Consider the map, called the boundary-bulk map,
\be\label{iota-def-eq}
\tau^E:\Hom_{\CC}(E,E)\rTo{\sim}\Tr_{\CC}(E^{\vee}\ot E)\rTo{\Tr_{\CC}(c_E)}
\Tr_{\CC}(\De)=HH_*(\CC),
\ee
where the first arrow is the isomorphism of Lemma \ref{Tr-lem}(ii).
Then $\ch(E)=\tau^E(\id_E)$.
\end{prop}

\Pf . By definition, $\ch(E)$ is obtained by considering the map on Hochschild homology
induced by the functor $\one_E:\per_{dg}(k)\to\per_{dg}(\CC)$. Thus, in our construction of $(\one_E)_*$
we should consider the canonical maps $k=\De_k\to E\ot_{\CC} E^\vee\simeq \Hom_{\CC}(E,E)$ 
and $E^\vee\ot_k E\to \De_{\CC}$.
The former corresponds to the identity element in $\Hom_{\CC}(E,E)$ while the latter is the map 
\eqref{cE-map}.
It follows that the map on Hochschild homology 
$$(\one_E)_*:k=HH_*(k)\to HH_*(\CC)$$
sends $1$ to $\tau^E(\id_E)$.
\ed

We will need the following compatibility of the maps \eqref{cE-map} with dg-functors.

\begin{lem}\label{cE-adjunction-lem} 
Let $F=T_X:\per_{dg}(\CC)\to\per_{dg}(\DD)$ be  a dg-functor, where $X\in\per_{dg}(\CC^{op}\ot\DD)$, and let
$G=T_Y:\per_{dg}(\DD)\to \per_{dg}(\CC)$ the right quasi-adjoint functor to $F$ given by \eqref{G-fun-eq}.

\noindent
(i) For $A\in\per_{dg}(\CC)$ we have a natural isomorphism 
\begin{equation}\label{A-F-G-dual-isom}
A^\vee\circ G=Y\ot_{\CC}A^\vee\rTo{\sim} (A\ot_{\CC} X)^\vee=F(A)^\vee
\end{equation}
in $\per(\DD^{op})$.

\noindent
(ii) For $A\in\per_{dg}(\CC)$ we have
the following  commutative diagram in the category $\per(\CC^{op}\ot\DD)$ whose objects are
viewed as functors $\per_{dg}(\CC)\to\per_{dg}(\DD)$
\begin{diagram}
F\circ (A^\vee\ot A) &\simeq A^\vee\ot F(A)\rTo{\eps}&\bigl(F(A)^\vee\ot F(A)\bigr)\circ F\\
&\rdTo{F\circ c_A}&\dTo{c_{F(A)}\circ F}\\
&& F
\end{diagram}
where $\eps$ is given by the composition
$$A^\vee\ot F(A)\to \bigl(A^\vee\ot F(A)\bigr)\circ G\circ F\rTo{\eps'\circ F} \bigl(F(A)^\vee\ot F(A)\bigr)\circ F$$
in which 
$$\eps':\bigl(A^\vee\ot F(A)\bigr)\circ G\rTo{\sim} F(A)^\vee\ot F(A)$$
 is the isomorphism induced by \eqref{A-F-G-dual-isom}.
\end{lem}

\Pf . (i) This isomorphism is given as the following composition
\begin{align*}
&(X^T\ot_{\CC}A^\vee)(D)
\rTo{\sim}\Hom_{\CC^{op}}(A,X^T(D,?))\rTo{\sim}
\Hom_{\CC^{op}}(A,\Hom_{\DD^{op}}(X,h_D))\simeq\\
&\Hom_{\DD^{op}}(A\ot_{\CC}X,h_D)\simeq (A\ot_{\CC}X)^\vee(D).
\end{align*}

\noindent
(ii) Let us denote by $F_A$ the functor associated with the kernel
$A^\vee\ot F(A)\in\per(\CC^\vee\ot\DD)$:
\begin{equation}\label{F-A-fun-eq}
F_A=\Hom_\CC(A,?)\ot F(A):\per_{dg}(\CC)\to \per_{dg}(\DD).
\end{equation}
Consider the diagram
\begin{diagram}
F\circ[A^\vee\ot A]&\simeq F_A\rTo{}& F_A\circ G\circ F&\rTo{\eps'\circ F} & [F(A)^\vee\ot F(A)]\circ F\\
\dTo{F\circ c_A}&&\dTo{\phi\circ G\circ F}&&\dTo{c_{F(A)}\circ F}\\
F &\rTo{}& F\circ G\circ F&\rTo{}& F 
\end{diagram}
where $\phi:F_A\to F$ is the morphism given by the natural maps
\begin{equation}\label{phi-F-A-eq}
\Hom_\CC(A,?)\ot F(A)\to\Hom_\DD(F(A),F(?))\ot F(A)\to F,
\end{equation}
and the unmarked horizontal arrows are induced by the quasi-adjunction maps $\De_\CC\to G\circ F$ and
$F\circ G\to\De_\DD$ (see \eqref{quasi-adj-maps}). The left square is commutative because $F\circ c_A$
corresponds to $\phi$ under the isomorphism $F\circ[A^\vee\ot A]\simeq F_A$.
The commutativity of the right square reduces to the commutativity of the
diagram
\begin{diagram}
\Hom_\CC(A,GF(?))\ot F(A)&\rTo{}&\Hom_\DD(F(A), FGF(?))\ot F(A)&\rTo{}&\Hom_\DD(F(A),F(?))\ot F(A)\\
                                 &\rdTo{}&\dTo{}&&\dTo{}\\
                                 &&FGF&\rTo{}&F
                                 \end{diagram}
in which the composition of the arrows in the first line is equal to $\eps'\circ F$.
\ed

\subsection{Generalized abstract \HRR\ Theorem}
\label{cardy-sec}

It is known that a Calabi-Yau dg-category gives rise to an open-closed 2d TQFT
(see \cite{Costello}). One of the equations of the open-closed TQFT, the so called
Cardy condition, can be viewed as a generalization of the \HRR\ formula. It was observed
in \cite[Thm.\ 15]{CalW} that the Cardy condition can be stated without assuming the Calabi-Yau
property (in \cite{CalW} this condition is called the "Baggy Cardy Condition").
Here we prove a categorical version of the Cardy condition for an arbitrary dg-category $\CC$ such that
$\per_{dg}(\CC)$ is saturated.

\begin{thm}\label{cardy-prop}
For a pair of objects $A,B\in\per_{dg}(\CC)$ 
and elements $\a\in\Hom_{\CC}(A,A)$, $\b\in\Hom_{\CC}(B,B)$
we have
\begin{equation}\label{cardy-cat-eq}
\lan \tau^{A^{\vee}}(\a^{\vee}),\tau^B(\b)\ran=\sTr_k(m_{\a,\b}),
\end{equation}
where $\a^{\vee}\in\Hom_{\CC^{op}}(A^{\vee},A^{\vee})$ is induced by $\a$, and
$m_{\a,\b}$ is the endomorphism
\begin{equation}\label{m-a-b-eq}
m_{\a,\b}:\Hom_{\CC}(A,B)\to\Hom_{\CC}(A,B):f\mapsto (-1)^{|\a|\cdot |\b|+|\a|\cdot|f|}\cdot \b\circ f\circ \a.
\end{equation}
\end{thm}

Note that the \HRR\ formula \eqref{HRR-cat-eq} is obtained by setting $\a=\id_A$, $\b=\id_B$
and using Proposition \ref{ch-lem}. In the case when $\per(\CC)=D^b(X)$ is the bounded
derived category of
coherent sheaves on a smooth proper scheme $X$ over $k$, 
the formula \eqref{cardy-cat-eq} along with the expression of $\tau^{\FF}$ in
terms of the Atiyah class of $\FF\in D^b(X)$ can be found in \cite{Ram}.

The proof of \eqref{cardy-cat-eq}
will be based on the following result on compatibility of the boundary-bulk maps $\tau^A$ with the functoriality
\eqref{F*-eq} of Hochschild homology.

\begin{lem}\label{tau-funct-lem} 
Let $F:\per_{dg}(\CC)\to\per_{dg}(\DD)$ be a dg-functor.
Then for an object $A\in\per_{dg}(\CC)$ the following diagram is commutative
\begin{diagram}
\Hom_{\CC}(A,A) &\rTo{F}&\Hom_{\DD}(F(A),F(A))\\
\dTo{\tau^A}&&\dTo{\tau^{F(A)}}\\
HH_*(\CC)&\rTo{F_*}&HH_*(\DD)
\end{diagram}
\end{lem}

\Pf . Recall that the definition of $F_*$ uses the natural
isomorphism $\Tr_{\CC}(G\circ F)\simeq\Tr_{\DD}(F\circ G)$ (see Lemma \ref{Tr-lem}(i)),
where $G:\per_{dg}(\DD)\to\per_{dg}(\CC)$ is the right quasi-adjoint dg-functor to $F$.
Consider the dg-functor $F_A:\per_{dg}(\CC)\to\per_{dg}(\DD)$ given by \eqref{F-A-fun-eq}
and the morphism of dg-functors $\phi:F_A\to F$ given by \eqref{phi-F-A-eq}.
We have a commutative diagram
\begin{equation}\label{FA-F-G-diagram}
\begin{diagram}
\Tr_{\CC}(G\circ F_A)&\simeq&\Tr_{\DD}(F_A\circ G)\\
\dTo{\Tr_{\CC}(\id_G\circ\phi)}&&\dTo_{\Tr_{\DD}(\phi\circ\id_G)}\\
\Tr_{\CC}(G\circ F)&\simeq&\Tr_{\DD}(F\circ G)
\end{diagram}
\end{equation}
Let $\psi:A^{\vee}\ot A\to G\circ F_A$ be the morphism induced by the quasi-adjunction
$$\Hom_{\CC}(A,?)\ot A\to\Hom_{\CC}(A,?)\ot GF(A)\simeq G\circ F_A(?)$$
and $\eta:F_A\circ G\to F(A)^{\vee}\ot F(A)$ the morphism
induced by the quasi-adjunction
$$F_A\circ G(?)\simeq\Hom_{\CC}(A,G(?))\ot F(A)\to\Hom_{\DD}(F(A),?)\ot F(A).$$
These morphisms fit into the following commutative diagrams in $\per(\CC^{op}\ot\CC)$ and $\per(\DD^{op}\ot\DD)$,
respectively:
\begin{diagram}
A^{\vee}\ot A&\rTo{\psi}&G\circ F_A\\
\dTo{c_A}&&\dTo{\id_G\circ\phi}\\
\De_{\CC}&\rTo{}&G\circ F,
\end{diagram}
\begin{diagram}
F_A\circ G&\rTo{\eta}& F(A)^{\vee}\ot F(A)\\
\dTo{\phi\circ\id_G}&&\dTo{c_{F(A)}}\\
F\circ G&\rTo{}&\De_{\DD}
\end{diagram}
Now we apply $\Tr_{\CC}$ and $\Tr_{\DD}$ to these diagrams and use \eqref{FA-F-G-diagram}
to get the result.
\ed

\noindent
{\it Proof of Theorem \ref{cardy-prop}.} 
We apply Lemma \ref{tau-funct-lem} to the functor
$$\Tr^{dg}_{\CC}:\per_{dg}(\CC^{op}\ot\CC)\to\per_{dg}(k),$$
the object $A^{\vee}\ot B\in \per(\CC^{op}\ot\CC)$, and an element
$\a^{\vee}\ot \b\in \Hom_{\CC^{op}\ot\CC}(A^{\vee}\ot B,A^{\vee}\ot B)$.
We have $\Tr^{dg}_{\CC}(A^{\vee}\ot B)=\Hom_{\CC}(A,B)$ and 
$$\Tr^{dg}_{\CC}(\a^{\vee}\ot \b)=m_{\a,\b}$$
(the sign in \eqref{m-a-b-eq}
appears from the definition of $\Tr$, see Remark \ref{Tr-sign-rem}). 
It remains to use the fact that the maps $\tau^A$ are compatible with 
the K\"unneth isomorphism:
$$\tau^{A^{\vee}\ot B}=\tau^{A^{\vee}}\ot\tau^B.$$
\ed

\section{Matrix factorizations}
\label{sec:mf}
In this section we collect some facts
about categories of \mf s (see e.g. \cite{Dyck} and \cite{KhR} for more information).

\subsection{Categories of \mf s}

Let $R$ be a commutative algebra over a field $k$.
We fix an element $\w \in R$ which will be called
the \emph{potential}. 

\begin{defn}
A \emph{\mf} of the potential $\w$ over $R$ is a 
pair 
\be\label{mf-defn}
(E,\d_E)=(E^\zer \lrarrow{{\d_\zer}}{{\d_\one}} E^\one),
\ee
 where
\begin{itemize}
\item  $E=E^\zer\oplus E^\one$ is a $\Zt$-graded finitely generated projective $R$-module, and
\item $\d_E\in \End^\one_R(E)$ is an odd (i.e.\ of degree $\one \in \Zt$)
endomorphism of $E$, such that $\d_E^2=\w\cdot \id_E$.  
\end{itemize}
Even though $\d_E^2\ne 0$, we will still call $\d_E$ a ``differential''.
\end{defn}

If $E^\zer$ and $E^\one$ are free $R$-modules with chosen bases, 
the differential $\d_E$ can be represented by
a block matrix 
\be
 D= \left( \begin{array}{cc}
0 & D^\one\\ D^\zer & 0
 \end{array} \right), 
\ee
such that the matrices $D^\zer$ and $D^\one$ give a factorization of 
the potential: 
$$
 D^\zer D^\one = D^\one D^\zer = \w\cdot\mathbf{I}.
$$

To a potential $\w \in R$ we associate a $\Zt$-dg-category $\MF(\w)=\MF(R,\w)$
whose objects are \mf s of $\w$ over $R$. 
The morphisms from $\bar{E}=(E,\d_E)$ to  $\bar{F}=(F,\d_F)$
are elements of the $\Z/2$-graded module of 
$R$-linear homomorphisms
$$
\Homb_{\w}(\bar{E},\bar{F}):=\Hom_{\Mod_R}(E,F)
=\Hom_{\Zt-\Mod_R}(E,F)\oplus \Hom_{\Zt-\Mod_R}(E,F[\one]).
$$  
The $\Zt$-graded dg-structure on  $\MF(R,\w)$
is given by the differential $d$ defined on $f\in\Homb_{\w}(\bar{E},\bar{F})$ 
as 
\begin{equation}\label{diff-mf-eq}
d f = \d_F\circ f - (-1)^{|f|} f \circ \d_E~.
\end{equation}
For $\bar{E}, \bar{F}\in\MF(\w)$ we set
$$\Hom_\w(\bar{E},\bar{F})=H^*(\Homb_\w(\bar{E},\bar{F}),d).$$

Let
$$\MFt(R,\w)=H^0\MF(R,\w)$$
denote the homotopy category
associated with the dg-category $\MF(R,\w)$. By definition, 
morphisms in this category are chain maps up to homotopy, i.e,
the spaces $\Hom^0_\w(\bar{E},\bar{F})$.
The homotopy category $\MFt(R,\w)$ of \mf s is naturally
triangulated (see e.g.~\cite{Orlov})
with the shift functor induced from the functor
$T:\MF(R,\w)\to \MF(R,\w)$ given by
$$
T(E,\d_E)=(E[1],-\d_E), \quad T(f)=f[1], \mathrm{\ for \ }
f\in \Homb_{\w}(\bar{E},\bar{F}).
$$
The category $\MFt(R,\w)$ can be viewed as a full triangulated subcategory
of the perfect derived category of \mf s $\per(\MF(R,w))$ and is isomorphic to
it when $R$ is regular and complete by \cite[Thm.\ 5.7]{Dyck} (see \cite{O-compl} for
a more general statement). 

For a pair of elements $a,b \in R$
denote by $\{ a, b \}$ the \mf\ of the potential $\w=ab$ given by
\be
 \{ a, b \}=(R \lrarrow{a}{b} R)\in \MF(R,ab).
\ee

An easy computation shows that 
\begin{equation}\label{ab-hom}
\begin{array}{l}
\Hom^0_{\w}(\{a,b\},\{a,b\})\simeq R/(Ra+Rb),\\
\Hom^1_{\w}(\{a,b\},\{a,b\})\simeq\{(x,y)\in R\oplus R\ |\ ax=by\}/R\cdot (b,a),
\end{array}
\end{equation}
provided $\w$ is not a zero divisor in $R$.
In particular, the identity map of  $\{ a, b \}$ is homotopic to zero 
(and so $\{ a, b \}$ represents the zero object in $\MFt(R,ab)$)
if and only if the ideal $(a,b)$ coincides with $R$. Also, if the pair
$(a,b)$ forms a regular sequence then $\Hom^1_{\w}(\{a,b\},\{a,b\})=0$.

If $G$ is a finite group of automorphisms of $R$ which fixes the potential
$\w$, one defines the $G$-equivariant $\Zt$-graded dg-category of \mf s
$\MF_G(\w)$ (and the corresponding homotopy category $\MFt_G(\w)$), by
requiring that all modules and morphisms should be $G$-equivariant (see e.g., \cite{Ve}).
In other words, in \eqref{mf-defn}
$E$ should be a $\Zt$-graded finitely generated projective $R$-module $E$ equipped with
a compatible $G$-action, and $\de_E$ has to be $G$-equivariant. Morphisms between
$G$-equivariant \mf s $\bar{E}$ and $\bar{F}$ should also be compatible with the action of $G$, so 
\begin{equation}\label{G-hom-eq}
\Homb_{\MF_G(\w)}(\bar{E},\bar{F})=\Homb_{\w}(\bar{E},\bar{F})^G.
\end{equation}

\subsection{Tensor product and duality}
\label{sec:kernels}

The \emph{tensor product} $\bar{E}\otimes_R \bar{E}'$ of 
two \mf s $\bar{E}=(E,\d_E)\in \MF(R,\w)$  and  $\bar{E}'=(E',\d_{E'}) \in \MF(R,\w')$
is defined as the pair
\be
(E\otimes_R E', \d_E\otimes \id_{E'} + J_E\otimes \d_{E'}),
\ee
where
$J_E=(-1)^{|\cdot|}$ is the grading operator. Since $\de_E\ot\id_{E'}$ and $J_E\ot\de_{E'}$
anticommute, the tensor product $\bar{E}\otimes_R \bar{E'}$ is a \mf\ of
the potential $\w+\w'$.

Note that we have a natural commutativity isomorphism in $\MF(\w+\w')$:
\begin{equation}\label{comm-isom-eq}
\bar{E}\ot_R\bar{E}'\simeq \bar{E}'\ot_R\bar{E}:e\ot e'\mapsto (-1)^{|e|\cdot |e'|}e'\ot e.
\end{equation}

The following definition was introduced in~\cite{BGS} (see also~\cite{KhR}).
\begin{defn}
Let $\mathbf{a}=(a_1,\ldots,a_n)$ and $\mathbf{b}=(b_1,\ldots,b_n)$ be two
$n$-tuples of elements of $R$. The \mf\
\be
\label{eq:kz_mf}
\{\mathbf{a},\mathbf{b}\} := \{a_1,b_1\}\otimes_R\ldots \otimes_R
 \{a_n,b_n\}
\ee
of the potential
$$ \w = \mathbf{a}\boldsymbol{\cdot}\mathbf{b}=a_1 b_1+\ldots+a_n b_n
$$
is called \emph{the Koszul \mf\ } corresponding to the pair
$(\mathbf{a},\mathbf{b})$.
\end{defn}

More explicitly,
$\{\mathbf{a},\mathbf{b}\}$ is
isomorphic as a $\Zt$-graded $R$-module
to the Koszul complex 
\be\label{kos-com}
K_\bullet=\left({\bigwedge}_R^\bullet(R^{n}),\d\right),
\ee
where the differential is given by 
$$\d=(\sum_{j=1}^na_je_j)\we ?+\iota(\sum_{j=1}^n b_je_j^*).$$ 
Here $(e_j)$ is the standard basis of $R^n$, $(e_j^*)$ is the dual basis of
the dual $R$-module, $\iota$ denotes the contraction operator.
The $\Z/2$-grading is induced by the $\Z$-grading of $K_\bullet$. 

We also need a $G$-equivariant version of Koszul matrix factorizations. Suppose $G$ is a finite group
acting on $R$ and $\w\in R$ is $G$-invariant. Assume we have
a finite-dimensional $G$-module
$V$ and a pair of $G$-invariant elements $\phi\in V\ot R$ and $\psi\in V^*\ot R$ 
such that 
$$\lan\psi,\phi\ran=\w.$$
Then we get a structure of a
$G$-equivariant \mf\ of $\w$ on the corresponding Koszul complex
\be\label{G-kos-com}
K_{\bullet}(V)={\bigwedge}_R^{\bullet}(V\ot R),
\ee
where the differential is given by
$$\de=\phi\we ?+\iota(\psi).$$
We denote this $G$-equivariant \mf\ by $\{\phi,\psi\}$. Choosing a basis $(e_i)$ in $V$
we can write $\phi=\sum_i e_i\ot a_i$, $\psi=\sum_i e_i^*\ot b_i$,
where $(e_i^*)$ is the dual basis in $V^*$. Then after forgetting a $G$-action we obtain
$$\{\phi,\psi\}=\{\mathbf{a},\mathbf{b}\}.$$

For a pair of potentials $\w\in R$ and $\w'\in R'$ 
we have a natural external tensor product functor 
\be\label{ext-tensor-eq}
\MF(R,\w) \otimes \MF(R',\w') \to \MF(R\otimes R',\w\oplus\w'),
\ee
where $\w\oplus\w' := \w\otimes 1_{R'} + 1_R\otimes\w' \in R\otimes R'$,
sending a pair $\bar{E}\in\MF(R,\w)$, $\bar{E'}\in\MF(R,\w')$ to
$$\bar{E}\boxtimes\bar{E'}=(E\ot_k E', \de_E\ot\id+J_E\ot\de_{E'}).$$ 


To a \mf\ 
$\bar{E}=(E^\zer \lrarrow{{\d_\zer}}{{\d_\one}} E^\one)$
of the potential $\w\in R$ we associate the following {\it dual \mf} of $-\w$ 
\be\label{dual-mf-def}
\bar{E}^*=((E^\zer)^* \lrarrow{{\d_\one^*}}{{-\d_\zer^*}} (E^\one)^*)\in
\MF(R,-\w),
\ee
where for a projective $R$-module $P$ we set $P^*=\Hom_R(P,R)$.
In other words, for $e^*\in E^*$ and $e\in E$ we have
$$\lan \de_{E^*}(e^*), e\ran=(-1)^{|e^*|}\lan e^*,\de_E(e)\ran
$$
which is the usual sign rule for defining the adjoint operator in the $\Z/2$-graded context
(since $|\de_E|=1$).
The functor $\bar{E}\mapsto\bar{E}^*$ gives an equivalence 
\be\label{dual-op-eq}
\MF(R,\w)^{op}\simeq\MF(R,-\w).
\ee

For any $\bar{E},\bar{F} \in \MF(R,\w)$
the tensor product $\bar{F}\otimes_R \bar{E}^*$ 
is a $\Z/2$-graded complex (since it is a \mf\ of the potential $\w+(-\w)=0$)
and we have a natural isomorphism of complexes of $R$-modules
\be\label{hom-dual-tens-eq}
\bar{F}\otimes_R \bar{E}^* \wt{\to} \Homb_\w(\bar{E},\bar{F})
\ee
(note that our choice of sign in the definition \eqref{dual-mf-def} is compatible with this isomorphism).

For $\bar{E}\in\MF(R,\w)$, $\bar{F}\in\MF(R,\w')$ we have an isomorphism
$$(\bar{E}\ot\bar{F})^*\wt{\to}\bar{F}^*\ot\bar{E}^*$$
of \mf s of $-\w-\w'$, given by the natural pairing between the $R$-modules $F^*\ot_R E^*$ and $E\ot_R F$:
$$\lan f^*\ot e^*, e\ot f\ran=e^*(e)\cdot f^*(f),$$
where $e\in E$, $f\in F$, $e^*\in E^*$, $f^*\in F^*$.
Note that the double dual is 
$$(\bar{E}^*)^*=\bar{E}_-:=(E,-\de_E),$$ 
so the natural isomorphism
$\bar{E}\to(\bar{E}^*)^*$ is given by the grading operator $J_E$.
This is compatible with the sign convention
$$\lan e,e^*\ran=(-1)^{|e|}\lan e^*,e\ran$$
for $e\in E$, $e^*\in E^*$ (the pairing is nonzero only if $|e|=|e^*|$).

Thus, we see that for $\bar{E},\bar{F}\in\MF(R,\w)$ there is
a natural pairing
\be\label{tensor-pairing-eq}
\Homb_{\w}(\bar{E},\bar{F})\ot \Homb_{\w}(\bar{F},\bar{E})\to R: \ A\ot B\mapsto \sTr_R(A\circ B)
\ee
inducing the perfect duality between these complexes,
where $\sTr_R$ denotes the {\it supertrace} of an $R$-linear operator $C$ (it is 
equal to $\tr_R(C_{00})-\tr_R(C_{11})$, where $C_{ii}:E_i\to E_i$ are the components of $C$).
Indeed, if we use the identifications 
$\Homb_{\w}(\bar{E},\bar{F})\simeq \bar{F}\ot_R\bar{E}^*$ and
$\Homb_{\w}(\bar{F},\bar{E})\simeq \bar{E}\ot_R\bar{F}^*$ then the
above pairing corresponds to the standard pairing
$$\lan f\ot e^*, e\ot f^*\ran=(-1)^{|f|}\lan e^*,e\ran\cdot \lan f^*,f\ran.
$$
In the case when $R$ is a regular local $k$-algebra of dimension $n$, and $\w$ is an isolated
singularity (see sec. \ref{sec:generator} below),
by Grothendieck duality, the above duality induces an isomorphism
$$\Hom_{\w}(\bar{E},\bar{F})^*\simeq\Hom_{\w}(\bar{F},\bar{E})\ot_R \om_{R/k}[n]$$
which means that $\bar{E}\mapsto \bar{E}\ot_R\om_{R/k}[n]$ is a Serre functor on
$\MFt(R,\w)$ (see~\cite{Aus}, \cite{Bu}). 
An explicit formula for the duality trace maps is 
proved in \cite{Murfet}.


\subsection{Matrix factorizations and modules over
hypersurface singularities}
\label{sec:mcm}

From now on we assume that $R$ is a regular local $k$-algebra with
the maximal ideal $\m$ and the residue field $R/\m=k$.
Given a minimal set of generators $x_1,\ldots,x_n$ of $\m$, we denote the
corresponding derivations of $R$ by $\del_i$.
We will be mostly interested in the case when
$R$ is the ring of formal power series 
 $k[[x_1,\ldots,x_n]]$.

Matrix factorizations of a potential $\w\in \m\subset R$
naturally arise in the study of maximal \CM\ modules
over the hypersurface algebra $S=R/\w$.
(A \emph{maximal  \CM\ module} over a commutative Noetherian local 
ring is a module whose depth is equal to the Krull dimension of the ring.)  
If 
$$
(E,\d_E)=(E^\zer \lrarrow{{\d_\zer}}{{\d_\one}} E^\one)
$$
is a \mf\  of $\w\in R$ then, since
$\w\cdot E \subset \im(\d_E) $,
the $R$-module $\Cok(\d_\one)$ is naturally an $S$-module.
Eisenbud showed in~\cite{Eis} that $\Cok(\d_\one)$ is 
a maximal \CM\ module
and that any maximal \CM\ module over $S$ can be obtained this way. 
Moreover, he proved that the functor
$\mathrm{Coker}: \MF(R,\w) \to \Mod_S$ induces an equivalence 
of categories 
\be
\label{eq:eisen_eq}
 \MFt(R,\w) \to \MCMs(S),
\ee
where $\MCMs(S)$ is the \emph{stable category of maximal \CM\ $S$-modules}
(the quotient of the full subcategory of maximal \CM\ $S$-modules
$\MCM(S)\subset \Mod_S$ modulo free $S$-modules).

Buchweitz~\cite{Bu} extended this result to a characterization of the
\emph{stabilized derived category} of $S$,
$$
\Dbs(S)= \Db(S)/\Dbp(S),
$$
where
$\Db(S)$ is the bounded derived category of all complexes of $S$-modules with 
finitely generated cohomology
and
$\Dbp(S)$ is 
the full triangulated subcategory of $\Db(S)$
of perfect complexes (i.e.\ complexes quasi-isomorphic 
to a bounded complex of free $S$-modules).
He proved that the natural functor
$\MCMs(S)\to \Dbs(S)$ is an equivalence of categories which, 
combined with \eqref{eq:eisen_eq}, 
induces an equivalence of triangulated categories
\be
\label{eq:eq_buch}
\MFt(R,\w) \to \Dbs(S)
\ee
(see \cite{Orlov} for a generalization).
Thus, to every finitely generated $S$-module $M$ there corresponds
a \mf\ $M^{\st}$ in $\MFt(\w)$ such that its image under~(\ref{eq:eq_buch})
is isomorphic to $M$ in $\Dbs(S)$.
Following Dyckerhoff~\cite{Dyck}, 
we call $M^\st$ the \emph{stabilization} of $M$.
(This term reflects the fact proved in~\cite{Eis} 
that $M$ has a free resolution which is stably 2-periodic.)

For a large class of $S$-modules the stabilizations are provided by 
Koszul matrix factorizations~(\ref{eq:kz_mf}).

\begin{prop}\label{koszul-mf-prop}
[\cite[Sec.\ 7]{Eis}, cf.\ also~\cite[Cor.\ 2.7]{Dyck}]
Let $I$ be an ideal in $R$ generated by 
a regular sequence $\mathbf{b}=(b_1,\ldots,b_n)$, and
let $\mathbf{a}=(a_1,\ldots,a_n)$
be an $n$-tuple of elements of $R$ such that $\w=\mathbf{a}\cdot\mathbf{b}$.
Then the Koszul matrix factorization 
$\{\mathbf{a},\mathbf{b}\}=(\bigwedge_R^\bullet(R^n),\de)\in \MFt(R,\w)$
gives the stabilization of the $S$-module $R/I$. More precisely,
the natural map of $S$-modules
\begin{equation}\label{cokers-map}
\coker(\de_1)\to \coker\left({\bigwedge}^1 \rTo{\iota(\sum_j b_je_j^*)}{}
{\bigwedge}^0\right)=R/I
\end{equation}
induced by the projection $\bigwedge^\bullet\to\bigwedge^0$ becomes an isomorphism in $\Dbs(S)$.
\end{prop}

\Pf . Let us set $K^i=\bigwedge_R^i(R^n)$ and let 
$$\de_i(a):K^i\to K^{i+1}, \ \de_i(b):K^i\to K^{i-1}$$
denote the components of the differential $\de$. 
We can view $\coker(\de_1:K^{odd}\to K^{even})$ as the $0$th cohomology of the
total complex of the bicomplex $L^{\bullet,\bullet}$ concentrated on two diagonals
$i+j=0$ and $i+j=-1$ given by
$$L^{-i,i}=K^{2i} \ \text{ and } L^{-i,i-1}=K^{2i-1}$$
with the differentials $\de_{2i-1}(b):L^{-i,i-1}\to L^{-i,i}$ and
$\de_{2i-1}(a):L^{-i,i-1}\to L^{-i+1,i-1}$. Consider the spectral sequence
of this bicomplex that starts by taking the cohomology in the horizontal direction.
The $E_1$-term is given by
$$E_1^{-i,i}\simeq \coker(\de_{2i+1}(b)),$$
$$E_1^{-i,i-1}\simeq \ker(\de_{2i-1}(b))\simeq\coker(\de_{2i+1}(b)) \text{ for } i\ge 1,$$
where the last isomorphism uses the exactness of the Koszul complex associated with 
the regular sequence $(b_1,\ldots,b_n)$. The differential 
$E_1^{-i,i-1}\to E_1^{-i,i}$ for $i\ge 1$ is given by the multiplication with $\w$
(recall that both modules can be identified with $\coker(\de_{2i+1}(b))$).
But $\coker(\de_{2i+1}(b))\simeq\ker(\de_{2i-1}(b))$ is isomorphic to a submodule of a free $R$-module,
so the multiplication by $\w$ on it is injective. It follows that the $E_2$-term is concentrated on
the diagonal $i+j=0$ and
$$E_2^{-i,i}\simeq \coker(\de_{2i+1}(b)\ot_R S) \text{ for } i\ge 1,$$
while $E_2^{0,0}\simeq R/I$. It remains to observe that the complex
$K^\bullet\ot_R S$ with the differential $\de_\bullet(b)$ computes $\Tor_\bullet^R(S,k)$,
so it is exact in the terms $K^i$ with $i\ge 2$, since the projective dimension of $S$ as $R$-module is
equal to $1$. Hence, for $i\ge 1$ the $S$-module
$E_2^{-i,i}$ admits a finite free $S$-resolution, so it becomes zero in $\Dbs(S)$.
\ed

\subsection{Generators and Hochschild homology for \mf s}
\label{sec:generator}

From now on we assume that $R$ is the ring of formal power series 
$k[[x_1,\ldots,x_n]]$. 

\begin{defn}
An element $\w\in \m$ is called an \emph{isolated singularity} if
its \emph{Tyurina algebra} $R/(\w,\del_1 \w,\ldots,\del_n \w)$
has finite dimension over $k$. 
\end{defn}

It is well known (see \cite[Prop.\ (1.2)]{Loo}) that in this case the
\emph{Milnor ring} 
$$\AA_{\w}=R/(\del_1\w,\ldots,\del_n\w)$$ 
is also finite dimensional (this uses the assumption that $k$ has characteristic zero).

If $\w$ is an isolated singularity then it is finitely determined, i.e., $\w$ is determined by
its finite jet up to a change of variables 
(see \cite{Mather}, or in algebraic setting \cite{Hir} and \cite{CS}). In particular, there exists an
automorphism $\phi$ of $R$ over $k$ such that $\phi(\w)$ is a polynomial in $x_1,\ldots,x_n$.
This simple fact will allow us to extend the results of \cite[Sec.\ 6]{Dyck}, proved for
algebraic potentials, to the case when $R=k[[x_1,\ldots,x_n]]$.

Let $\MF^{\infty}(R,\w)$ be the $\Zt$-dg-category 
of matrix factorizations involving free $R$-modules
of possibly infinite rank and let
$\HMF^{\infty}(R,\w)=H^0\MF^\infty(R,\w)$ be the corresponding homotopy category.
The Yoneda embedding $\HMF(R,\w)\to D(\MF(R,\w))$ extends to a functor
\begin{equation}\label{ioneda-functor-eq}
\HMF^{\infty}(R,\w)\to D(\MF(R,\w)): \bar{E}\mapsto h_{\bar{E}}:\bar{F}\mapsto\Homb_\w(\bar{F},\bar{E})).\end{equation}
 
\begin{thm}\label{gen-thm}
Let $\w\in R$ be an isolated singularity. Let us view
$R/\m\simeq k$ as a module
over $S=R/\w$, and let $k^\st\in\MF(R,\w)$ be the stabilization of $k$.
Then the Yoneda functor
\eqref{ioneda-functor-eq} is an equivalence
and the category $D(\MF(R,\w))$
is compactly generated by $k^{\st}$. The $\Zt$-dg-category $\per_{dg}(\MF(R,\w))$ is saturated.
 \end{thm}

\Pf . Let $R_0$ be the localization of $k[x_1,\ldots,x_n]$ at the maximal ideal corresponding to the origin.
In the case when $R$ is replaced by $R_0$
the theorem follows from \cite[Thm.\ 4.1, Thm.\ 5.2, Prop.\ 6.3]{Dyck}.
By finite determinacy of $\w$ we can assume that $\w\in R_0\sub R$.
But $R$ is the completion of $R_0$, so we have quasi-equivalences
$$\per_{dg}(\MF(R_0,\w))\to\per_{dg}(\MF(R,\w)), \  \MF^{\infty}(R_0,\w)\to \MF^{\infty}(R,\w)$$
(the first equivalence follows from \cite[Lem.\ 5.6, Thm.\ 5.7]{Dyck}, and the second is proved
in the same way using \cite[Thm.\ 5.2]{Dyck}).
\ed
 
As a consequence of the above theorem, one gets a dg Morita equivalence of the category  
$\MF(R,\w)$ with the $\Zt$-dg-algebra 
$$A=\Homb^*_{\w}(k^{\st},k^{\st}),$$ 
so the computation of
$HH_*(\MF(R,\w))$ reduces to that of $HH_*(A)$. However, technically
it is more convenient to use the category of \mf s for the doubled potential
\be\label{tilde-w-eq}
\wt{\w}:=\w(y_1,\ldots,y_n)-\w(x_1,\ldots,x_n),
\ee
which is an element of the ring 
$$R^e:=R\hat{\ot}_k R=k[[x_1,\ldots,x_n,y_1,\ldots,y_n]].$$

The arguments of \cite[Sec.\ 6.1]{Dyck} imply that the external tensor product functor 
\eqref{ext-tensor-eq}
induces an equivalence 
\begin{equation}\label{double-equiv}
\MFt(R^e,\wt{\w})\simeq\per(\MF(R,\w)^{op}\ot\MF(R,\w)).
\end{equation}
Indeed, by the finite determinacy of $\w$ we can assume that $\w\in R_0\sub R$, where
$R_0$ is the localization of $k[x_1,\ldots,x_n]$ at the maximal ideal corresponding to the origin.
Then we get from \cite[Sec.\ 6.1]{Dyck} an
equivalence of derived categories
$$D(\MF(R_0\ot_k R_0,\wt{\w}))\simeq D(\MF(R_0,\w)^{op}\ot\MF(R_0,\w)).$$
Since $R$ is the completion of $R_0$ and $R^e$ is the completion of $R_0\ot_k R_0$,
by passing to perfect subcategories we obtain \eqref{double-equiv}.

Note that the equivalence \eqref{double-equiv} sends a \mf\ $K$ of $\wt{\w}$ to the 
$\MF(\w)-\MF(\w)$-bimodule
\begin{equation}\label{hK-eq}
h_K(\bar{E},\bar{F})=
\Homb_{\wt{\w}}(\bar{E}^*\boxtimes\bar{F}, K)\simeq (\bar{E}\boxtimes\bar{F}^*)\ot_{R^e} K,
\end{equation}
where we denote by $\bar{E}^*\boxtimes\bar{F}\in\MF(R^e,\wt{\w})$ 
the completed external tensor product of $\bar{E}^*$ and
$\bar{F}$.

Let us show, following Dyckerhoff \cite[Sec.\ 6]{Dyck}, how the above results lead
to the computation of the Hochschild homology of $\MF(R,\w)$ for
an isolated singularity $\w$. To apply the definition \eqref{hochschild-hom-def}
we need an explicit object in $\MF(R^e,\wt{\w})$ 
representing the diagonal $\MF(R,\w)-\MF(R,\w)$-bimodule. 
As shown in
\cite[Prop.\ 6.3]{Dyck} (see also \cite[Prop.\ 23]{KhR}), the stabilized diagonal 
$\De^{\st}\in\MF(R^e,\wt{\w})$ is such an object.
Explicitly, this is the Koszul \mf\ (see \eqref{eq:kz_mf}, \eqref{kos-com})
\be\label{De-st-eq}
\De^{\st}=\{\De_1\w,\ldots,\De_n\w;y_1-x_1,\ldots,y_n-x_n\}=(K_{\bullet}^{\De},\de_K),
\ee
associated with the decomposition
$$\wt{\w}=\w(y)-\w(x)=\sum_{j=1}^n\De_j\w\cdot (y_j-x_j),$$
where 
\begin{equation}\label{Wj-eq}
\De_j\w=\frac{\w(x_1,\ldots,x_{j-1},y_j,y_{j+1}\ldots,y_n)-\w(x_1,\ldots,x_{j-1},x_j,y_{j+1},\ldots,y_n)}
{y_j-x_j}\in R^e.
\end{equation}
For later calculations it is important to 
note that the difference derivative 
$\De_j\w$ does not depend on $y_1,\ldots,y_{j-1}$ and on $x_{j+1},\ldots,x_n$ and that
\be\label{De-restr-eq}
\De_j\w|_{y=x}=\pa_j \w.
\ee 

Proposition \ref{koszul-mf-prop} implies that 
the explicit isomorphism of the $\MF(\w)-\MF(\w)$-bimodule 
$h_{\De^{\st}}$ with the diagonal bimodule $\De_{\MF(\w)}$
is given by the map
\begin{equation}\label{ga-map-eq}
h_{\De^{\st}}(\bar{E},\bar{F})\simeq (\bar{E}\ot\bar{F}^*)\ot_{R\ot R}K_\bullet^{\De}\rTo{\ga} 
\bar{E}\ot_R \bar{F}^*\simeq\Homb_\w(\bar{F},\bar{E})=\De_{\MF(\w)}(\bar{E},\bar{F}).
\end{equation}
Here $\ga$ is induced by the composition 
$$K_{\bullet}^{\De}\to K_0^{\De}=R^e\to R^e/\JJ_{\De}=R,$$
where $\JJ_{\De}=(y_1-x_1,\ldots,y_n-x_n)\sub R^e$ is the ideal of the diagonal.

On the other hand, we claim that the composition
$$\HMF(R^e,\wt{\w})\to \per(\MF(\w)^{op}\ot\MF(\w))\rTo{\Tr}\per(k),$$
where $\Tr=\Tr_{\MF(\w)}$ is the functor \eqref{tr-per-eq},
is isomorphic to the functor of restriction to the diagonal
$$K\mapsto K|_{y=x}=K\ot_{R^e}R.$$
Indeed, Theorem \ref{gen-thm} implies that $\HMF(R^e,\wt{\w})$ is generated by objects of
the form $\bar{F}^*\boxtimes\bar{E}$, where $\bar{E},\bar{F}\in\MF(\w)$,
so our claim follows from the isomorphism
$$\Tr(\bar{F}^*\ot\bar{E})\simeq\Homb_\w(\bar{F},\bar{E})\simeq
(\bar{F}^*\boxtimes\bar{E})|_{y=x}.$$

Hence, the complex calculating
$HH_*(\MF(R,\w))$ is obtained by restricting $(K^{\De}_\bullet,\d_K)$ to the diagonal $y=x$.
This gives the usual Koszul differential for the regular sequence $\pa_1\w,\ldots,\pa_n\w$,
so one obtains an isomorphism 
$$\ga_{\bx}:HH_*(\MF(R,\w))\rTo{\sim} \AA_{\w}[n]$$
(see \cite[Thm.\ 6.6]{Dyck}). Multiplying it with the top-degree form $d\bx=dx_1\we\ldots\we dx_n$
we obtain a {\it canonical} isomorphism
\be\label{H-w-eq}
\ga=\ga_{\bx}\ot d\bx: HH_*(\MF(R,\w)) \rTo{\sim} H(\w):=\AA_{\w}\ot \om_{R/k}[n],
\ee
where $\om_{R/k}=\Om^n_{R/k}$. 
To see that this isomorphism does not depend on a choice of a system of parameters $\bx=(x_1,\ldots,x_n)$ we note
that the restriction of $K_{\bx}=(K_\bullet^{\De},\de_K)$ to the diagonal can be naturally identified with the
complex $(\Om^\bullet_{R/k}, d\w\we ?)$. Any other 
regular system of parameters $\bx'=(x'_1,\ldots,x'_n)$ for $R$ gives rise to a different
Koszul \mf\ $K_{\bx'}$ of $\wt{\w}$ representing the identity functor.
There is a homotopy equivalence between $K_{\bx}$ and $K_{\bx'}$ which is compatible with the identification
of the restrictions of $K_{\bx}$ and $K_{\bx'}$ to the diagonal with $(\Om^\bullet_{R/k}, d\w\we ?)$.
This implies that the isomorphism \eqref{H-w-eq} does not depend on a choice of a system of parameters.
In particular, it is unchanged when we permute the
variables and is compatible with the action of the group of symmetries of $\w$.

\subsection{$G$-equivariant Hochschild homology}
\label{sec:G-hoch}

Here we present an equivariant version of the results of section \ref{sec:generator}. 
Let $G$ be a finite group acting on the
algebra $R=k[[x_1,\ldots,x_n]]$ by automorphisms (identical on $k$), 
and let $\w\in R$ be a $G$-invariant isolated singularity. Recall that we denote $S=R/(\w)$.
Let $R\#G$ (resp., $S\#G$) denote the twisted group ring of $G$ over $R$ (resp., $S$).

There is a $G$-equivariant version of the equivalence 
\eqref{eq:eq_buch},
$$\MFt_G(R,\w)\simeq\Dbs_G(S),$$
where on the right one takes the quotient of the bounded
derived category of finitely generated $S\#G$-modules
by the subcategory generated by $S\#G$-modules that are free over $S$ (see \cite{Ve}).
In particular, to every $G$-equivariant
finitely generated $S$-module $M$ there corresponds naturally a $G$-equivariant \mf\ $M^{\st}$,
the {\it stabilization} of $M$.

Let us describe explicitly the stabilization of $k=R/\m$.
Since we work in characteristic zero, the surjective map of $G$-modules
$$\m\to\m/\m^2=:V$$
admits a $G$-equivariant splitting $s:V\to\m$. 
Let $\psi\in (V^*\ot R)^G$ be the element corresponding to $s:V\to R$.
The $G$-equivariant map
$$\lan ?,\psi\ran:V\ot R\to \m$$
is surjective. Hence, we can find a $G$-invariant element
$\phi\in V\ot R$ such that $\lan\phi,\psi\ran=\w$.
Then $\{\phi,\psi\}$ is a $G$-equivariant Koszul \mf\ (see \eqref{G-kos-com}), and
\be\label{G-k-st}
k^{\st}=\{\phi,\psi\}.
\ee

It will be convenient for what follows to choose a 
basis $(e_1,\ldots,e_n)$ in the space $V=\m/\m^2$ and use the 
elements 
\be\label{x-e-eq}
x_i=s(e_i)\in\m, \ i=1,\ldots,n,
\ee
as a new set of variables (recall that
we work with formal power series). 
With respect to these coordinates
$G$ acts on $R$ by linear transformation, which will be 
frequently used in the rest of the paper.

Let $\MF_G^{\infty}(\w)$ 
denote the $\Zt$-dg-category of $G$-equivariant \mf s of free $R$-modules
of possibly infinite rank.

 \begin{thm}\label{G-gen-thm}
The category $\MF_G^{\infty}(\w)$ is quasi-equivalent to
the $\Zt$-dg-derived category of $\MF_G(R,\w)$ and is compactly
generated by the $G$-equivariant \mf\ $k^\st\otimes k[G]$. 
 \end{thm}

\Pf .  For a \mf\ $\bar{E}\in\MF_G^{\infty}(\w)$
we have
$$\Hom_{\w}^*(k^{\st}\ot k[G],\bar{E})^G\simeq\Hom_{\w}^*(k^{\st},\bar{E}).$$
Since $k^{\st}$ is a generator of $\MFt^{\infty}(\w)$, this immediately implies that
$k^{\st}\ot k[G]$ is a generator of the homotopy category $\MFt^{\infty}_G(\w)$.
Now the same argument as in the proof of Theorem 5.2 of \cite{Dyck} implies that $\MFt^{\infty}_G(\w)$
is quasi-equivalent to the derived category of $\MF_G(\w)$.
\ed

Since the finite determinacy also holds in the equivariant setting (see \cite[Prop.\ 5.1]{Roberts} or 
\cite[Lem.\ 1.2]{Wall}), there exists a $G$-equivariant automorphism $\phi$ of $R$ over $k$ such that
$\phi(\w)$ is a polynomial .
Thus, as in the non-equivariant case the external tensor product functor induces an equivalence
of $\per(\MF_G(R,\w)^{op}\ot\MF_G(R,\w))$ with $\HMF_{G\times G}(R^e,\wt{\w})$. Hence,
to proceed with the calculation of Hochschild homology we should work
with kernels which are
$G\times G$-equivariant matrix factorizations of $\wt{\w}$.

 To obtain an analog of the diagonal factorization \eqref{De-st-eq} we start with an element
\be\label{psi-eq}
\psi_{\De}=\sum_{i=1}^n e_i^*\ot (y_i-x_i)\in V^*\ot R^e,
\ee
and set 
$$\phi_{\De}=\sum_{i=1}^n e_i\ot \De_i\w\in V\ot R^e,$$ 
so that $\lan\phi_{\De},\psi_{\De}\ran=\wt{\w}=\w(y)-\w(x)$.
If we view $R^e$ as a $G$-module via the diagonal action on $(x_1,\ldots,x_n,y_1,\ldots,y_n)$
then $\psi_{\De}$ will be $G$-invariant (by our choice of variables \eqref{x-e-eq}).
Replacing $\phi_{\De}$ with 
\be\label{G-phi-eq}
\phi^G_{\De}=|G|^{-1}\cdot \sum_{g\in G}g\cdot \phi_{\De}=:\sum_{i=1}^n e_i\ot \wt{\De}_i\w
\ee
we obtain a $G$-equivariant Koszul \mf 
\be\label{G-De-st-eq}
\De_G^{\st}:=\{\phi^G_{\De},\psi_{\De}\}=\{\wt{\De}_1\w,\ldots,\wt{\De}_n\w;y_1-x_1,\ldots,y_n-x_n\}
\ee
(see \eqref{G-kos-com}) which is isomorphic to $\De^{\st}$ after forgetting the $G$-equivariant structure.
Recall that the map $V\ot R^e\to V\ot R$
(restriction to the diagonal $y=x$)
 sends $\phi_{\De}$ to the $G$-invariant element $\sum_{i=1}^n e_i\ot \pa_i\w$
(see \eqref{De-restr-eq}). Hence, we still have
\be\label{G-De-restr-eq}
\wt{\De}_i\w|_{y=x}=\pa_i\w.
\ee
Now we define a $G\times G$-equivariant \mf\ of $\wt{\w}$ 
\be\label{GG-De-st-eq}
\De_{G\times G}^{\st}:=\bigoplus_g (\id\times g)^*\De_G^{\st}.
\ee


\begin{prop} The $G\times G$-equivariant \mf\ $\De_{G\times G}^{\st}$
of $\wt{\w}$ is a kernel for the identity functor on $\MFt_G^{\infty}(R,\w)$.
Together with Theorem \ref{G-gen-thm}
this implies that the $\Zt$-dg-category $\per_{dg}(\MF_G(R,\w))$ is saturated.
\end{prop}

\Pf . The proof is parallel to the non-equivariant case (see \cite[Prop.\ 6.3]{Dyck}).
\ed

Now in order to calculate $HH_*(\MF_G(R,\w))$ we have to compute $\Tr^G(\De_{G\times G}^{\st})$, 
where 
$$\Tr^G:\MF_{G\times G}(R^e,\wt{\w})\to\per(k)$$ 
is the categorical trace functor \eqref{Tr-def-eq}.
It follows from the description \eqref{G-hom-eq} 
of morphisms in $\MF_G(R,\w)$ as $G$-invariants  
that $\Tr^G$ is given by the restriction to the diagonal $y=x$ followed by taking $G$-invariants.
Thus, we need first to compute the restriction of each summand $(\id\times g)^*\De_G^{\st}$ to the diagonal.
We use an obvious isomorphism
\be\label{Ga-g-restr}
(\id\times g)^*\De_G^{\st}|_{\De}\simeq\De_G^{\st}|_{\Ga_g},
\ee
where 
$$\Ga_g=\{(x,y)\ |\ y=gx\}\sub\Spec(R\hat{\ot}_k R)$$ 
is the graph of the action of $g$ on $\Spec(R)$.

Furthermore, for a given element $g\in G$ we can choose variables 
$(x_1,\ldots,x_n)$ in such a way that 
\begin{equation}\label{g-eq}
g(x_1,\ldots,x_n)=(\ell_1,\ldots,\ell_k,x_{k+1},\ldots,x_n),
\end{equation}
where $\ell_1,\ldots,\ell_k$ are linear forms in $x_1,\ldots,x_k$, and 
$\Span(x_{k+1},\ldots,x_n)$ is exactly the subspace of $g$-invariants in
$\Span(x_1,\ldots,x_n)$, i.e., the forms $(\ell_1-x_1,\ldots,\ell_k-x_k)$ are linearly independent.
Let us consider the restriction of $\w$ to the subspace of $g$-invariants: 
$$\w_g(x_{k+1},\ldots,x_n):=\w|_{x_1=\ldots=x_k=0}.$$
In other words, we take the image of $\w$ under the projection 
$R\to R^g=R/(x_1,\ldots,x_k)$.

\begin{lem}\label{Ga-g-lem} 
(i) For $j=k+1,\ldots,n$
we have
\be\label{De-diag-g-eq}
\wt{\De}_j\w|_{\Ga_g\cap\De}=\pa_j\w_g,
\ee 
and the sequence $(\pa_{k+1}\w_g,\ldots,\pa_n\w_g)$
in $R^g$ is regular, i.e., $\w_g\in R^g$ is an isolated singularity.

\noindent
(ii) The cohomology of the complex $\De_G^{\st}|_{\Ga_g}$ maps isomorphically to
$$H(\w_g)=\AA_{\w_g}\cdot dx_{k+1}\we\ldots\we dx_n[n-k],$$
via the restriction to $R/(x_1,\ldots,x_k)\ot e_{k+1}\we\ldots \we e_n$ and passing to the quotient modulo
$(\pa_{k+1}\w_g,\ldots,\pa_n\w_g)$ (where $dx_j$ gets identified with $e_j$).
\end{lem}

\Pf .
(i) The equation \eqref{De-diag-g-eq} 
follows by setting $x_1=\ldots=x_k=0$ in the identity $\wt{\De}_j\w|_{\De}=\pa_j\w$ for $j>k$ (see \eqref{G-De-restr-eq}).
To check that $\w_g$ is an isolated singularity let us differentiate the equation $\w(gx)=\w(x)$
with respect to the variables $x_1,\ldots,x_k$ and substitute $x_1=\ldots=x_k=0$. Using
linear independence of the forms $\ell_1-x_1,\ldots,\ell_k-x_k$ we derive
that 
$$\pa_i \w|_{x_1=\ldots=x_k=0}=0 \text{  for }i=1,\ldots,k.$$
It follows that any critical point of $\w_g$ is also a critical point of $\w$ which implies our claim.

\noindent
(ii) 
We have 
$\De_G^{\st}=A_{\bullet}\ot B_{\bullet}$, where
$A_{\bullet}$ (resp., $B_{\bullet}$)
is the matrix factorization of $\wt{\w}_A=\sum_{i=1}^k (y_i-x_i)\wt{\De}_i\w$
(resp., $\wt{\w}_B=\sum_{j=k+1}^n (y_j-x_j)\wt{\De}_j\w$),
associated with this decomposition. Since $\wt{\w}|_{\Ga_g}=\w(gx)-\w(x)=0$ and
$\wt{\w}_B|_{\Ga_g}=0$, it follows that $\wt{\w}_A|_{\Ga_g}=0$. Hence, we have a decomposition
into the product of two complexes 
$$\De_G^{\st}|_{\Ga_g}=A_{\bullet}|_{\Ga_g}\ot B_{\bullet}|_{\Ga_g},$$
where $A_{\bullet}|_{\Ga_g}$ is the complex associated with the decomposition
$$0=\w_A|_{\Ga_g}=\sum_{i=1}^k(\ell_i-x_i)\cdot \wt{\De}_i\w,$$ 
and $B_{\bullet}$ is isomorphic to 
the Koszul complex for the sequence $(\wt{\De}_{k+1}\w|_{\Ga_g},\ldots,\wt{\De}_n\w|_{\Ga_g})$, shifted by 
$n-k$ (and twisted by $dx_{k+1}\we\ldots\we dx_n$). Since the forms $\ell_1-x_1,\ldots,\ell_k-x_k$
are linearly independent, it follows that the 
cohomology of $A_{\bullet}$ is isomorphic to $k$ and is concentrated in the term $A_0$, so
$$H^*(\De_G^{\st}|_{\Ga_g})\simeq H^*(B_{\bullet}|_{\De\cap\Ga_g}).$$
By part (i), this reduces to the Koszul complex for the regular
sequence $(\pa_{k+1}\w_g,\ldots,\pa_n\w_g)$.
\ed

The above computation leads to the following decomposition of the equivariant Hochschild homology
of $\MF_G(R,\w)$, similar to the decompositions of Hochschild and cyclic 
(co)homology of twisted group algebras
(see \cite{Burg, FT, Lorenz}).

\begin{thm}\label{G-hoch-thm} Let $\w\in R$ be an isolated singularity, invariant under a finite
group of automorphisms $G$. 
Then 
we have
\begin{equation}\label{G-hoch-eq}
HH_*(\MF_G(R,\w))\simeq (\bigoplus_{g\in G} H(\w_g))^G,
\end{equation}
where an element $g\in G$ acts on $\bigoplus_{g\in G} H(\w_g)$ by sending
$H(\w_h)$ to $H(\w_{ghg^{-1}})$.
\end{thm}

We will also need the following result.

\begin{lem}\label{choices-rem}
The isomorphism \eqref{G-hoch-eq} does not depend on a choice of a $G$-invariant 
element $\phi^G_{\De}\in V\ot R^e$ such that 
\begin{equation}\label{De-phi-psi-eq}
\wt{\w}=\lan \phi^G_{\De},\psi_{\De}\ran,
\end{equation}
where $\psi_{\De}$ is given by \eqref{psi-eq}.
\end{lem}

\Pf .
It is enough to check that for two different choices $\phi$ and $\phi'$ of the element
$\phi^G_{\De}$ satisfying \eqref{De-phi-psi-eq}
there exists
an isomorphism  $\{\phi,\psi_{\De}\}\simeq\{\phi',\psi_{\De}\}$
in $\HMF(\wt{\w})$ such that
its restriction to $\De\cap\Ga_g$ induces the identity map on cohomology under the identification
of Lemma \ref{Ga-g-lem}(ii). Indeed, the equivalence \eqref{eq:eq_buch} implies that  
the restriction of any \mf\ $(E,\de)$ of $\wt{\w}$
to $\De\cap\Ga_g$ depends functorially only on the module $\coker(\de_1)$ 
(cf. \cite[Lem.\ 4.2]{Dyck}). 
We claim
that for any pair of Koszul \mf s $\{\mathbf{a},\mathbf{b}\}$ and $\{\mathbf{a'},\mathbf{b}\}$ of $\w\in R$ with
regular $\mathbf{b}=(b_1,\ldots,b_n)$ 
there exists an isomorphism $\a:\{\mathbf{a},\mathbf{b}\}\to\{\mathbf{a'},\mathbf{b}\}$ in $\HMF(\w)$
that induces the identity on $R/(b_1,\ldots,b_n)$.
This will imply the statement of the lemma because of the isomorphism 
\eqref{cokers-map}. To prove our claim we observe that because of the acyclicity of the Koszul
complex of the regular sequence $\mathbf{b}$
we can find an element $x=\sum_{i<j}x_{ij}e_i\we e_j\in \bigwedge^2_R(R^n)$ such that
$$\iota(\sum_{i=1}^n b_i e_i^*)(x)=\sum_{i=1}^n a'_ie_i - \sum_{i=1}^n a_ie_i,$$
where $\mathbf{a}=(a_1,\ldots,a_n)$ and $\mathbf{a'}=(a'_1,\ldots,a'_n)$. Now
the operator $\a=\exp(-x)\we ?$ provides the required isomorphism.
\ed

\begin{rem} E.~Segal in \cite{Segal} obtains a similar result for 
a slightly modified version of Hochschild homology of $\MF_G(R,\w)$ (he used a
version of the standard complex,
where direct sums are replaced
by direct products). It is not clear a priori how to identify his version with 
Hochschild homology defined in the usual way.
\end{rem}

\section{Chern character and boundary-bulk maps}
\label{Chern-sec}

Let $\bar{E}=(E,\de_E)$ be a matrix factorization of $\w$.
Our goal in this section
is to compute the image of the Chern character $\ch(\bar{E})\in HH_*(\MF(R,{\w}))$
(see \eqref{ch-def-eq}) under the isomorphism $HH_*(\MF(R,\w))\simeq H(\w)$ (see \eqref{H-w-eq}).
More generally, we will compute explicitly the boundary-bulk map \eqref{iota-def-eq}
$$\tau^{\bar{E}}:\Hom_{\w}(\bar{E},\bar{E})\to HH_*(\MF(R,\w))\simeq H(\w). $$

\subsection{Technical lemmas}
\label{sec:complex}
Let $\De^{\st}$ be the diagonal matrix factorization of $\wt{\w}=\w(y)-\w(x)$
(see \eqref{De-st-eq}). We will reduce the problem of computing the map $\tau^{\bar{E}}$
to finding a special element in the $\Zt$-graded complex of $R^e$-modules
\begin{equation}\label{L-complex}
L_{\bullet}:=\De^{\st}\ot_{R^e} p_2^*\bar{E}^*\otimes_{R^e} p_1^*\bar{E}
\simeq\Homb_{\wt{\w}}(\bar{E}^*\boxtimes\bar{E},\De^{\st}),
\end{equation}
where the last isomorphism is the combination of the isomorphism \eqref{hom-dual-tens-eq} and
the duality \eqref{tensor-pairing-eq}. 
Recall that $\De^{\st}=(K_{\bullet},\de_K)$, where
$$K_{\bullet}=\bigoplus_{i=0}^n K_i, \text{ and } K_i={\bigwedge}^i_{R^e}((R^e)^n).$$ 
By \eqref{H-w-eq}, $H(\w)$ is isomorphic to the cohomology of the complex
$\De^{\st}|_{y=x}=(K_{\bullet}|_{y=x},\de_K|_{y=x})$ concentrated in the term $K_n|_{y=x}$.
Let us denote by 
$$\pi:K_n|_{y=x}\to H(\w)$$ 
the corresponding projection.

\begin{lem}\label{str-lem} 
Let $D=\sum_{j=0}^n D_j\in L_{even}$ 
be a closed even element in the complex \eqref{L-complex} such that 
$$D_0|_{y=x}=1\ot\id_E\in K_0|_{y=x}\ot E^*\ot_R E,$$
where $D_i\in K_i\ot_R p_2^*E^*\ot_R p_1^*E$.
Then the value of the boundary-bulk map
on $\a\in\Hom_{\w}(\bar{E},\bar{E})$ is given by
$$\tau^{\bar{E}}(\a)=\pi(\sTr(D_n|_{y=x}\circ\a)).$$
In particular, for the Chern character of $\bar{E}$ we have
$$\ch(\bar{E})=\pi(\sTr(D_n|_{y=x})).$$
\end{lem}

\Pf . Recall that the external tensor product functor induces an equivalence 
$$\per(\MF(R,\w)^{op}\ot\MF(R,\w))\simeq \HMF(R^e,\wt{\w}),$$
such that the diagonal bimodule $\De_{\MF(\w)}$ corresponds to the diagonal \mf\ $\De^{\st}$
and the functor $\Tr:\per(\MF(R,\w)^{op}\ot\MF(R,\w))\to\per(k)$ corresponds to the
functor of restriction to the diagonal.
Since elements of the complex $L_{\bullet}$ are morphisms from $\bar{E}^*\boxtimes \bar{E}$
to $\De^{\st}=(K_\bullet,\de_K)$ in $\MF(R^e,\wt{\w})$, 
we have to find an element $c_{\bar{E}}^{\vee}\in H^0(L_{\bullet})$ inducing the canonical
map  \eqref{cE-map}
\begin{equation}\label{c-map-eq}
c_{\bar{E}}:\bar{E}^*\boxtimes \bar{E}\to \De^{\st}:\a\mapsto\lan c_{\bar{E}}^{\vee},\a\ran
\end{equation}
and then restrict it to the diagonal.
By definition, the map $c_{\bar{E}}$ is characterized by the functorial commutative diagram 
\begin{diagram}
(\bar{F}\ot_R\bar{E}^*)\ot(\bar{G}^*\ot_R\bar{E})\\
\dTo{\lan?,c_{\bar{E}}^{\vee}\ran}&\rdTo{t_E}\\
h_{\De^{\st}}(\bar{F},\bar{G})\simeq(\bar{F}\ot\bar{G}^*)\ot_{R\ot R}K_\bullet
&\rTo{\ga}&\bar{F}\ot_R\bar{G}^*
\end{diagram}
where $\ga$ is the map \eqref{ga-map-eq} and $t_E$ is induced by the evaluation map
$E^*\ot_R E\to R$.
Since $\ga$ is induced by the projection $K_{\bullet}\to K_0=R^e\to R^e/\JJ_{\De}=R$, 
we obtain that the $K_0$-component of $c_{\bar{E}}^{\vee}$ projects to $1\ot\id_E$ under the projection
$K_0=R^e\to R$.
We observe that the latter projection coincides with the composition of the natural maps
$$K_0\to \coker(\d_K:K_{odd}\to K_{even})\to \coker(K_1\stackrel{\d_K}{\to} K_0),$$
where the second arrow becomes an isomorphism in the stable category.
By the equivalence \eqref{eq:eisen_eq}, an element
in $\Hom_{\wt{\w}}^0(\bar{E}^*\boxtimes \bar{E},\De^{\st})$ is determined by
the induced map on cokernels. 
Hence, for a closed element
$D\in L_{\bullet}$ such that $D_0|_{y=x}=1\ot\id_E$,
the image of $D$ in $H^0(L_{\bullet})$ is equal to $c_{\bar{E}}$.
Thus, we can take $c_{\bar{E}}^{\vee}=D$ in \eqref{c-map-eq}. 
When computing the restriction to the diagonal we recall that 
$K_{\bullet}|_{y=x}$ can be identified with the shifted Koszul complex
for the regular sequence $\pa_1\w,\ldots,\pa_n\w$ (up to a twist by $dx_1\we\ldots\we dx_n$), 
so its cohomology is
concentrated in the $K_n|_{y=x}$ term. Hence, it remains to take the $K_n|_{y=x}$-component of
the induced map 
$$H^*(\bar{E}^*\ot \bar{E})\to H^*(K_{\bullet}|_{y=x}):\a\mapsto\lan D|_{y=x},\a\ran=\sTr(D|_{y=x}\circ\a).
$$ 
\ed

We will need the following result about Koszul complexes.

\begin{lem}\label{Koszul-lem}
Let $A$ be a commutative ring, $R=A[t_1,\ldots,t_n]$, and let 
$$K_{\bullet}(t_1,\ldots,t_n)=\left({\bigwedge}_R^{\bullet}(R^n), \de\right)$$
be the Koszul complex for the sequence $(t_1,\ldots,t_n)$. Here
$\de=\iota(\sum t_je_j^*)$, where $(e_1,\ldots,e_n)$ is the standard basis of $R^n$,
$(e_1^*,\ldots,e_n^*)$ is the dual basis of the module $(R^n)^*$. 
Consider the following $A$-submodule in $K_{\bullet}$
$$C=C(t_1,\ldots,t_n)=\sum_{i_1<\ldots<i_k, k\ge 1}R_{\ge i_1}\cdot e_{i_1}\we \ldots\we e_{i_k},$$
where $R_{\ge i}=A[t_i,t_{i+1},\ldots,t_n]\sub R$.
Then 
$$K_{\bullet}=\ker(\de)\oplus C.$$
\end{lem}

\Pf . Let us use induction in $n$. For $n=1$ we have $C=R\cdot e_1$ and $\ker(\d)=R$, so the statement is clear.
Suppose the statement holds for $n-1$. 
We have 
\begin{equation}\label{K-tensor-isom}
K_{\bullet}(t_1,\ldots,t_n)=K_{\bullet}(t_1)\ot_A K_{\bullet}(t_2,\ldots,t_{n}).
\end{equation}
Let $\de_1$ (resp., $\de_{2,\ldots,n}$) denote the differential in $K_{\bullet}(t_1)$
(resp., $K_{\bullet}(t_2,\ldots,t_{n})$).
Note that under the isomorphism \eqref{K-tensor-isom} we have a direct sum decomposition
$$C(t_1,\ldots,t_n)=\left(1\ot C(t_2,\ldots,t_n)\right)\oplus \left(C(t_1)\ot_A K_{\bullet}(t_2,\ldots,t_n)\right).$$
In other words, we can write every element of $K_{\bullet}(t_1,\ldots,t_n)$ in the form
\begin{equation}\label{koszul-element}
x=\sum_{i\ge 0}(t_1^i\ot a_i+(t_1^i\cdot e_1)\ot b_i),
\end{equation} 
where $a_i,b_i\in K_{\bullet}(t_2,\ldots,t_n)$, and this
element is in $C(t_1,\ldots,t_n)$ if and only if $a_i=0$ for $i>0$ and $a_0\in C(t_2,\ldots,t_n)$.
Note that for $f,g\in A[t_1]$ and $a,b\in K_{\bullet}(t_2,\ldots,t_n)$ we have
\begin{equation}\label{koszul-tensor-diff}
\de(f\ot a+(g\cdot e_1)\ot b)=f\ot\de_{2,\ldots,n}(a)+gt_1\ot b-(g\cdot e_1)\ot \de_{2,\ldots,n}(b).
\end{equation}
This implies that
$$gt_1\ot b= (g\cdot e_1)\ot\de_{2,\ldots,n}(b)\ \mod(\im\de),$$
Also, we have 
$$f\ot\ker(\de_{2,\ldots,n})\in\ker\de.$$
Thus, starting with an arbitrary element $x\in K_{\bullet}(t_1,\ldots,t_n)$ 
and applying the induction assumption to its decomposition
\eqref{koszul-element}
we can write $a_i=a'_i+c_i$, where $a'_i\in\ker(\de_{2,\ldots,n})$
and $c_i\in C(t_2,\ldots,t_n)$. Then we have
$$x= \sum_{i\ge 0} (t_1^i\ot c_i+(t_1^i\cdot e_1)\ot b_i)\ \mod(\ker\de).$$
Furthermore, for $i>0$ we have
$$t_1^i\ot c_i= (t_1^{i-1}\cdot e_1)\ot\de_{2,\ldots,n}(c_i)\ \mod(\ker\de),$$
which proves that $K_{\bullet}(t_1,\ldots,t_n)=\ker\de+C(t_1,\ldots,t_n)$.
On the other hand, if the element $x$ is in $C(t_1,\ldots,t_n)$, so that $a_i=0$ for $i>0$ and
$a_0\in C(t_2,\ldots,t_n)$, then the equation $\de x=0$ would 
give by \eqref{koszul-tensor-diff}
$$1\ot \de_{2,\ldots,n}(a_0)+\sum_{i\ge 0}t_1^{i+1}\ot b_i-\sum_{i\ge 0}(t_1^i\cdot e_1)\ot \de_{2,\ldots,n}(b_i)=0.$$
Since the first components of the three summands lie in complementary $A$-submodules we derive
that $\de_{2,\ldots,n}(a_0)=0$ and $b_i=0$ for all $i$. By induction assumption, this implies that
$a_0=0$ as well.
\ed

\subsection{Formula for the Chern character and the boundary-bulk map}
\label{sec:calc}

We keep the notation of section \ref{sec:complex}.
Now we are going to construct an element $D\in L_{even}$ satisfying the conditions of Lemma \ref{str-lem}.
Consider the operators
$$\d_{\De}=\iota(\sum_{j=1}^n (y_j-x_j)e_j^*) \text{ and }\ 
\d_{\w}=(\sum_{j=1}^n \De_j\w\cdot e_j)\we ?$$ 
on $K_{\bullet}^{\De}$ so that $\d_{K}=\d_{\De}+\d_{\w}$.
Let us fix an isomorphism 
\be\label{E-U-eq}
E\simeq U\ot R,
\ee
where $U$ is a $\Z/2$-graded vector space.
Then we can view the differential $\de_E$ of $\bar{E}$ as an $R$-valued odd endomorphism of $U$.
The decomposition \eqref{E-U-eq} induces an isomorphism
$$p_2^*\bar{E}^*\ot_{R^e} p_1^*\bar{E}\simeq U^*\ot U\ot R^e\simeq\End(U)\ot R^e$$
of \mf s of $-\wt{\w}=\w(x)-\w(y)$, where
the differential $\wt{\d}$ on $\End(U)\ot R^e$ acts by
$$\wt{\d}(M)=\de_E(x)\circ M-(-1)^{|M|}M\circ\de_E(y).$$
The complex $L_{\bullet}$ (see \eqref{L-complex}) can now be expressed as
$L_{\bullet}=K_{\bullet}^{\De}\ot\End(U)$, and its differential is given by
$$\d_L=(\d_{\De}+\d_{\w})\ot\id_{\End(U)}+J_K\ot\wt{\de},$$
where $J_K$ is the grading operator on the $\Zt$-graded space $K_\bullet$. 
Thus, if we write 
$$D=\sum_{j=0}^n D_j\in L_{even}$$ 
with $D_j\in K_j\ot \End(U)$ then the condition that
$D$ is $\d_L$-closed is equivalent to the system
\begin{equation}\label{D-system}
(\d_{\De}\ot\id)(D_{j+1})+(\d_{\w}\ot\id)(D_{j-1})+(-1)^j(\id\ot\wt{\de})(D_j)=0,
\end{equation}
for $j=0,\ldots, n$, where we set $D_{-1}=D_{n+1}=0$.
For brevity we will write $\d_{\De}$ (resp.,
$\d_{\w}$, resp., $\wt{\de}$) instead of $\d_{\De}\ot\id$ (resp., $\d_{\w}\ot\id$, resp., $\id\ot\wt{\de}$). 

\begin{lem}\label{D-sol-lem} 
There exists a solution $D$ of \eqref{D-system} with $D_0=1\ot\id_U\in R^e\ot\End(U)$.
Furthermore, let us consider the decomposition
$$D_j=\sum_{i_1<\ldots<i_j} e_{i_1}\we\ldots\we e_{i_j}\ot D_j(i_1,\ldots,i_j),$$ 
where $D_j(i_1,\ldots,i_j)\in R^e\ot\End(U)$. Then 
there exists a unique $D$  satisfying \eqref{D-system} such that
the coefficient $D_j(i_1,\ldots,i_j)$ does not depend on $y_k$ with $k<i_1$, for all $j$, $i_1<\ldots<i_j$.
\end{lem}

\Pf . We define $D_j$ inductively starting with $D_0=1\ot\id_U$.
If all $D_i$ for $i\le j$ are
already defined, then to find $D_{j+1}$ we have to solve
\begin{equation}\label{koszul-rec}
-\d_{\De}(D_{j+1})=\d_{\w}(D_{j-1})+(-1)^j\wt{\d}(D_j)
\end{equation}
obtained from \eqref{D-system}. Suppose first that $j=0$.
Then the equation becomes 
\begin{equation}\label{d-D1-eq}
\d_{\De}(D_1)=-1\ot\wt{\d}(\id_U)=1\ot (\de_E(y)-\de_E(x)).
\end{equation}
Since the right-hand side is zero for $y=x$, such $D_1$ exists. Applying
Lemma \ref{Koszul-lem}
to $A=k[x_1,\ldots,x_n]$ and $t_j=y_j-x_j$,
we can find a unique $D_1$ satisfying \eqref{d-D1-eq}
such that $D_1(i)$ does not depend on $y_j$ with $j<i$.
Now, assume that $j>0$. Then the argument is similar, but we have to check first that
the right-hand side of \eqref{koszul-rec} is $\d_{\De}$-closed. Indeed, using the same
equation for $j-1$ we get
\begin{equation}
\begin{array}{l}
\d_{\De}(\d_{\w}(D_{j-1})+(-1)^j\wt{\d}(D_j))=
\d_{\De}\d_{\w}(D_{j-1})+(-1)^j\wt{\d}\d_{\De}(D_j)= \nonumber\\
(\d_{\De}\d_{\w}+\wt{\d}^2)(D_{j-1})-(-1)^j\wt{\d}\d_{\w}(D_{j-2}).
\end{array}
\end{equation}
Using the identity $\d_{\De}\d_{\w}+\wt{\d}^2=-\d_{\w}\d_{\De}$ and applying \eqref{koszul-rec} for $j-2$
we can rewrite this as
$$
\d_{\w}[-\d_{\De}(D_{j-1})-
(-1)^j\wt{\d}(D_{j-2})]=(\d_{\w})^2(D_{j-3})=0
$$
as claimed.
\ed

\begin{lem}\label{restriction-lem} 
The unique solution $D$ of \eqref{D-system} constructed in Lemma \ref{D-sol-lem} satisfies 
\begin{equation}\label{D-recursion}
\begin{array}{l}
D_j(n-j+1,\ldots,n)|_{y_{n-j+1}=x_{n-j+1}}=\\
D_{j-1}(n-j+2,\ldots,n)\circ
\pa_{n-j+1}\d_E(x_1,\ldots,x_{n-j+1},y_{n-j+2},\ldots,y_n)
\end{array}
\end{equation}
for $j=1,\ldots,n$.
\end{lem}

\Pf . We use induction in $j$.
For $j=1$ the equation \eqref{D-recursion} takes form 
\begin{equation}\label{D-1}
D_1(n)|_{y_n=x_n}=\pa_n\de_E(x).
\end{equation}
To prove this let us rewrite \eqref{d-D1-eq} as
$$\sum_{j=1}^n(y_j-x_j)\cdot D_1(j)=\de_E(y)-\de_E(x).$$
Since $D_1(n)$ does not depend on $y_1,\ldots,y_{n-1}$, 
substituting $y_1=x_1,\ldots, y_{n-1}=x_{n-1}$ into this equation gives
$$(y_n-x_n)\cdot D_1(n)=\de_E(x_1,\ldots,x_{n-1},y_n)-\de_E(x_1,\ldots,x_{n}),$$
which immediately implies \eqref{D-1}.

Similarly, for $j>1$ 
substituting $y_i=x_i$ for $i=1,\ldots,n-j$ 
into the equation \eqref{koszul-rec} with $j$ replaced by $j-1$
and comparing the coefficients of $e_{n-j+2}\we\ldots\we e_n$ we get 
\be\label{D-j-eq}
\begin{array}{l}
(y_{n-j+1}-x_{n-j+1})D_j(n-j+1,\ldots,n)= [(-1)^j\wt{\d}(D_{j-1}(n-j+2,\ldots,n))-\\
\sum_{i=n-j+2}^j (-1)^{n-j+2-i}\cdot
\De_i\w\cdot D_{j-2}(n-j+2,\ldots,\hat{i},\ldots,n)]|_{y_1=x_1,\ldots,y_{n-j}=x_{n-j}}. 
\end{array}
\ee
Now recall that $\De_i\w$ and $D_m(i,\ldots)$ do not depend on $y_{n-j+1}$ for $i\ge n-j+2$.
Therefore, in the right-hand side of \eqref{D-j-eq}
only $\wt{\d}$ (that involves $\d_E(y)$) depends on $y_{n-j+1}$.
Hence, after differentiating the above equation with respect to $y_{n-j+1}$ and restricting to 
$y_{n-j+1}=x_{n-j+1}$ we obtain \eqref{D-recursion}.
\ed

\begin{thm} \label{ch-thm} 
Let $\w\in R=k[[x_1,\ldots,x_n]]$ be an isolated singularity, and let
$\JJ_{\w}=(\pa_1\w,\ldots,\pa_n\w)$. 
Then the boundary-bulk map on
an endomorphism $\a\in\Hom^*_{\w}(\bar{E},\bar{E})$ 
of a \mf\ $\bar{E}=(E,\d_E)\in \MF(R,\w)$ is equal to 
\be\label{tau-eq}
\tau^{\bar{E}}(\a)=\sTr_R(\pa_n \d_E\circ\ldots\circ\pa_1 \d_E\circ\a)\cdot dx_1\we\ldots\we dx_n \ \mod\ 
\JJ_{\w}\cdot dx_1\we\ldots\we dx_n,
\ee
where we view $\d_E$ and $\a$ as matrices with values in $R$
after choosing a basis in a free $R$-module $E$.

In particular, for the Chern character of $\bar{E}$ we have
\be\label{ch-eq2}
\ch(\bar{E})=\sTr_R(\pa_n\d_E\circ\ldots\circ\pa_{1}\d_E)\cdot
dx_{1}\we\ldots\we dx_n \ \mod\ \JJ_{\w}\cdot dx_1\we\ldots\we dx_n.
\ee
\end{thm}

\Pf .
Restricting \eqref{D-recursion} to the diagonal $y=x$ 
and combining the resulting equations for $j=1,\ldots,n$ we deduce that
$$D_n|_{y=x}=e_1\we\ldots\we e_n\ot\pa_n\d_E\circ \pa_{n-1}\d_E\circ\ldots\circ\pa_1\d_E.$$
Now the required formulas  follow 
from Lemma \ref{str-lem}.
\ed

\begin{cor}\label{perm-cor} With the notations of Theorem \ref{ch-thm}
the expression 
$$\sTr_R(\pa_n\d_E\circ\ldots\circ\pa_1\d_E\circ\a)
\cdot dx_1\we\ldots\we dx_n\ \mod\ \JJ_{\w}\cdot dx_1\we\ldots\we dx_n$$
is invariant under permutations of indices $(1,\ldots,n)$. Hence, we have
$$\tau^{\bar{E}}(\a)=(-1)^n\cdot\frac{1}{n!}\cdot
\sTr_R((d\d_E)^{\we n}\circ\a) \ \mod\ 
\JJ_{\w}\cdot dx_1\we\ldots\we dx_n.$$
\end{cor}

\Pf . The first assertion follows from the canonicity of the isomorphism \eqref{H-w-eq}.
The second formula is obtained upon expanding
$$(d\de_E)^{\we n}\circ\a=(\pa_1\de_E\cdot dx_1+\ldots+\pa_n\de_E\cdot dx_n)^{\we n}\circ\a$$
and using the first assertion. Note that the supertrace is zero unless $\a$ has degree $n \mod (2)$.
Thus, when we move the odd symbols $dx_i$ to the right, each swapping with $\a$ produces
the sign $(-1)^n$. Since we have $n$ such symbols, the resulting sign is $(-1)^{n^2}=(-1)^n$. 
\ed

\subsection{$G$-equivariant Chern character}

Now, we are going to discuss a $G$-equivariant version of the results of the previous section keeping the notation and assumptions of section \ref{sec:G-hoch}.

Let $\bar{E}=(E,\d_E)$ be a $G$-equivariant matrix factorization of $\w$.
We denote by 
$$\ch_G(\bar{E})\in HH_*(\MF_G(R,\w))$$ 
the
Chern character of $\bar{E}$ and by
$$\tau^{\bar{E}}_G:\Hom_{\w}(\bar{E},\bar{E})^G\to HH_*(\MF_G(R,\w))$$
the boundary-bulk map \eqref{iota-def-eq}.
Our goal is to compute explicitly for every $g\in G$ and
$\a\in \Hom_{\w}(\bar{E},\bar{E})^G$ the component 
$$\tau^{\bar{E}}(\a)_g\in H(\w_g)$$
of $\tau^{\bar{E}}_G(\a)\in HH_*(\MF_G(R,\w))$
with respect to the decomposition \eqref{G-hoch-eq}.

\begin{lem}\label{triv-lem} Let $E$ be an $R\#G$-module, free of finite rank as $R$-module.
There exists an isomorphism $E\simeq U\ot R$ of $R\#G$-modules, where
$U$ is a representation of $G$.
\end{lem}

\Pf . Let $U=E/\m E$. Since we work in characteristic zero, we can choose a $G$-equivariant
splitting $U\to E$ of the surjective map $E\to E/\m E=U$ of $G$-modules. 
The induced map $U\ot R\to E$ will be an isomorphism by the standard argument using
Nakayama Lemma.
\ed

We choose an isomorphism $E\simeq U\ot R$ as in the above Lemma and proceed
as in section \ref{sec:complex} 
to consider the complex 
$$L^G_{\bullet}=\De_G^{\st}\ot_{R^e}\End(U)$$ 
which is now
equipped with a $G$-action (diagonal on $R^e$). 
Note that Lemma \ref{str-lem} still holds with
the complex $L_{\bullet}$ replaced by $L^G_{\bullet}$ 
(which amounts to replacing the difference derivatives
$\De_i\w$ with their $G$-equivariant version $\wt{\De}_i\w$ defined
by \eqref{G-phi-eq}).

For an element $g\in G$ let us denote by 
$$\De_g:=\De\cap\Ga_g\sub\Spec(R^e)$$ 
the intersection
of the diagonal with the graph of $g$. If we choose variables $(x_1,\ldots,x_n)$ 
so that $g$ acts by the linear transformation
\eqref{g-eq} and $\Span(x_{k+1},\ldots,x_n)$ is exactly the subspace of $g$-invariants in
$\Span(x_1,\ldots,x_n)$, then $\De_g$ is given by equations
$$x_1=y_1=\ldots=x_k=y_k=0,\  x_{k+1}=y_{k+1},\ldots,x_n=y_n.$$

\begin{lem}\label{G-str-lem}
We have
\begin{equation}
\tau^{\bar{E}}(\a)_g=\sTr(D|_{\De_g}\circ g\circ\a)
\end{equation}
where $D$ defined as in Lemma \ref{str-lem}.
\end{lem}

\Pf . By Proposition \ref{ch-lem}, we have to compute the components of the canonical morphism
$$c^G_{\bar{E}}:\bar{E}^*\boxtimes\bar{E}\to\bigoplus_{g\in G}(\id\times g)^*\De_G^{\st}\simeq
\bigoplus_{g\in G}(g\times\id)^*\De_G^{\st}.$$
We claim that its component corresponding to $g=1$ coincides with the non-equivariant map 
\eqref{c-map-eq}. Indeed, let us consider the forgetful functor
$\Phi:\MF_G(\w)\to\MF(\w)$. Note that $\Phi$ is given by the kernel
$$\bigoplus_{g\in G}(g\times \id)^*\De_G^{\st}\in \MF_{G\times 1}(\wt{\w})\simeq
\per_{dg}(\MF_G(\w)^{op}\ot\MF(\w))$$ 
(where we do not use the action of $G$ on the second factor).
Applying Lemma \ref{cE-adjunction-lem} to the functor $F=\Phi$ and the object $A=\bar{E}$ 
we deduce that
$c^G_{\bar{E}}$, viewed as a morphism in $\HMF_{G\times 1}(\wt{\w})$, factors as the composition
$$\bar{E}^*\boxtimes\bar{E}\rTo{a_{\bar{E}^*}}
\bigoplus_{g\in G}g^*\bar{E}^*\boxtimes\bar{E}\rTo{\wt{c}_{\bar{E}}} \bigoplus_{g\in G}(g\times \id)^*\De_G^{\st},$$
where 
$a_{\bar{E}^*}$ is induced by the embedding of the component corresponding to $g=1$ in
$\bigoplus_{g\in G}g^*\bar{E}^*$, and 
$\wt{c}_{\bar{E}}$
is induced by $c_{\bar{E}}$. This immediately implies our claim.

Thus, via the isomorphism \eqref{Ga-g-restr} the restriction of the $g$-component of $c^G_{\bar{E}}$ to
the diagonal $y=x$ gets identified with the restriction of the non-equivariant map
$c_{\bar{E}}$ to the graph $\Ga_g$.
Hence, the $g$-component of $\tau^{\bar{E}}$ is obtained by
restricting $c_{\bar{E}}$ to $\Ga_g$ and using the isomorphism
$$\bar{E}^*\ot_R\bar{E}\wt{\to}\bar{E}^*\boxtimes\bar{E}|_{\Ga_g}$$
induced by the action of $g$ on $E$.
It remains to apply Lemma \ref{Ga-g-lem}.
\ed

Now we are ready to prove the formula for the $G$-equivariant Chern character and the boundary-bulk
map.

\begin{thm}\label{G-ch-thm} Fix an element $g\in G$. Choose variables $(x_1,\ldots,x_n)$ 
so that $g$ acts by the linear transformation
\eqref{g-eq} and $\Span(x_{k+1},\ldots,x_n)$ is exactly the subspace of $g$-invariants in
$\Span(x_1,\ldots,x_n)$.
Then for any $G$-equivariant \mf\ $\bar{E}=(E,\d_E)$ 
we have
\begin{equation}\label{G-tau-for}
\tau^{\bar{E}}_G(\a)_g=\sTr_{R^g}([\pa_n\d_E\circ\ldots\circ\pa_{k+1}\d_E\circ g\circ\a]|_{x_1=\ldots=x_k=0})|\cdot
dx_{k+1}\we\ldots\we dx_n \ \mod\ \JJ_{\w_g},
\end{equation}
where $\a\in\Hom_\w(\bar{E},\bar{E})^G$ and $R^g=R/(x_1,\ldots,x_k)$.

In particular,
\begin{equation}\label{G-ch-eq}
\ch_G(\bar{E})_g=\sTr_{R^g}([\pa_n\d_E\circ\ldots\circ\pa_{k+1}\d_E\circ g]|_{x_1=\ldots=x_k=0})\cdot
dx_{k+1}\we\ldots\we dx_n \ \mod\ \JJ_{\w_g}.
\end{equation}
\end{thm}

\Pf . First, note that the projection $HH_*(\MF_G(R,\w))\to H(\w_g)$ factors as the composition
$$HH_*(\MF_G(R,\w))\to HH_*(\MF_{G_0}(R,\w))\to H(\w_g),$$
where $G_0\sub G$ is the subgroup generated by $g$. Furthermore, when
computing the projection $HH_*(\MF_{G_0}(R,\w))\to H(\w_g)$ we can use
either $\De_G^{\st}$ or $\De_{G_0}^{\st}$ (by Lemma \ref{choices-rem} applied to $G_0$-equivariant
situation). Let us work with
$$\De_{G_0}^{\st}=\{ \wt{\De}_1\w,\ldots,\wt{\De}_n\w; y_1-x_1,\ldots,y_n-x_n\}$$
(see \eqref{G-De-st-eq}).
By Lemma \ref{G-str-lem}, the problem reduces to
proving the following identity:
$$D_{n-k}(k+1,\ldots,n)|_{\De_g}=
\pa_n\d_E(x)\circ\ldots\circ\pa_{k+1}\d_E(x)|_{x_1=\ldots=x_k=0}.$$
Since $G_0$ acting on $R$ changes only variables $x_1,\ldots,x_k$, it follows
that $\wt{\De}_j\w=\De_j\w$ for $j>k$. Now the above equation
can be verified by the same argument as in
Lemma \ref{restriction-lem}. 
\ed

We have also an equivariant version of Corollary \ref{perm-cor}.

\begin{cor}\label{G-perm-cor} With the notations of Theorem \ref{G-ch-thm}
the expression 
$$\sTr_{R^g}([\pa_n\d_E\circ\ldots\circ\pa_{k+1}\d_E\circ g\circ\a]|_{x_1=\ldots=x_k=0})
\cdot dx_{k+1}\we\ldots\we dx_n\ \mod\ \JJ_{\w}$$
is invariant under permutations of indices $(k+1,\ldots,n)$. Hence, we have
$$\tau^{\bar{E}}(\a)=(-1)^{n-k}\cdot \frac{1}{(n-k)!}\cdot
\sTr_{R^g}([(d\d_E)^{\we (n-k)}\circ g\circ\a]|_{x_1=\ldots=x_k=0}) \ \mod\ 
\JJ_{\w}\cdot dx_{k+1}\we\ldots\we dx_n.$$
\end{cor}

We have the following relation between the maps $\tau^{\bar{E}}$ and $\tau^{\bar{E}^*}$,
where $\bar{E}^*\in\MF_G(R,-\w)$ is the dual \mf\ to $\bar{E}$.

\begin{lem}\label{G-dual-lem}
Under the identification $H(-\w_g)=H(\w_{g^{-1}})$ 
for any $\bar{E}\in\MF_G(R,\w)$ 
we have
$$\tau^{\bar{E}^*}(\a^*)_g=\tau^{\bar{E}}(\a)_{g^{-1}},$$
where $\a\in\Hom_\w(\bar{E},\bar{E})^G$.
In particular,
$$\ch_G(\bar{E}^*)_g=\ch_G(\bar{E})_{g^{-1}}.$$
\end{lem}

\Pf .
Applying \eqref{G-tau-for} for $\bar{E}^*$ and the dual endomorphism $\a^*$ we get
\begin{equation}
\begin{array}{l}
\tau^{\bar{E}^*}(\a^*)_g=
\sTr(\pa_n \d^*_E\circ\ldots\circ\pa_{k+1} \d^*_E\circ (g^{-1})^*\circ \a^*)\cdot dx_{k+1}\we\ldots\we dx_n \ \mod\ 
\JJ_{\w}=\nonumber\\
(-1)^{{n-k\choose 2}+(n-k)|\a|}\sTr((\a\circ g^{-1}\circ\pa_{k+1}\d_E\circ\ldots\circ\pa_n\d_E)^*)\cdot 
dx_{k+1}\we\ldots\we dx_n
\ \mod\ \JJ_{\w}= \\
(-1)^{(n-k)|\a|}
\sTr(\a\circ g^{-1}\pa_{k+1}\d_E\circ\ldots\circ\pa_n\d_E)\cdot dx_n\we\ldots\we dx_{k+1}\ \mod\ \JJ_{\w}=\\
\sTr(\pa_{k+1}\d_E\circ\ldots\circ\pa_n\d_E\circ g^{-1}\circ \a)\cdot dx_n\we\ldots\we dx_{k+1}\ \mod\ \JJ_{\w},
\end{array}
\end{equation}
where we used the equalities $\sTr(M^*)=\sTr(M)$,
$\sTr(M_1\circ M_2)=(-1)^{|M|\cdot|N|}\sTr(M_2\circ M_1)$ and $\a\circ g^{-1}=g^{-1}\circ\a$. It remains to use Corollary \ref{G-perm-cor}
to see that this is equal to the right-hand side of \eqref{G-tau-for}.
\ed

\section{\HRR\ formula}

In this section we work out the explicit form of the categorical \HRR\ formula \eqref{HRR-cat-eq}
for the categories $\MF(R,\w)$ and $\MF_G(R,\w)$, where $\w\in R=k[[x_1,\ldots,x_n]]$ is an isolated singularity, using our calculation of the Chern characters from the previous section.

\label{form-sec}

\subsection{Non-equivariant case}\label{non-equiv-form-sec}

Recall that the Hochschild homology $HH_*(\MF(R,\w))$ can be identified with the $\Zt$-graded
space $H(\w)=\AA_{\w}\ot d\bx[n]$, where $\AA_{\w}$ is the Milnor ring of $\w$.
By \eqref{inv-De-eq},
the canonical bilinear form \eqref{pair-eq}
on 
\begin{equation}\label{HH-w-w-eq}
HH_*(\MF(R,\w))=H(\w)=H(-\w)=HH_*(\MF(R,-\w))
\end{equation}
is equal to the inverse of the tensor
$$\ch(\De^{\st})\in HH_*(\MF(R^e,\wt{\w}))\simeq HH_*(\MF(R,-\w))\ot  HH_*(\MF(R,\w)),
$$
where the diagonal \mf\ 
$\De^{\st}\in\MF(R^e,\wt{\w})$ is the kernel \eqref{De-st-eq}
representing the identity functor on $\MF(R,\w)$. 
We can calculate the Chern character of $\De^{\st}$
using the general formula \eqref{ch-eq2}
(this computation is contained implicitly in \cite[Sec.\ 5.1]{KR}).

For any $f\in k[[x_1,\ldots,x_n]]$ denote by 
$\De_jf\in k[[x_1,\ldots,x_n,y_1,\ldots,y_n]]$ the difference derivative 
\begin{equation}\label{De-f-eq}
\De_jf=\frac{f(x_1,\ldots,x_{j-1},y_j,y_{j+1}\ldots,y_n)-f(x_1,\ldots,x_{j-1},x_j,y_{j+1},\ldots,y_n)}
{y_j-x_j}\in R^e.
\end{equation}

\begin{prop}\label{form-prop} 
We have
$$\ch(\De^{\st})=(-1)^{{n\choose 2}}\cdot\det(\De_j(\pa_i\w))\in
\AA_{\wt{\w}}\cdot dx_1\we\ldots\we dx_n\we dy_1\we\ldots\we dy_n.$$
\end{prop}

\Pf . By \eqref{ch-eq2}, this reduces to equality (21) of \cite[Sec.\ 5.1]{KR}.
\ed

\begin{prop}\label{res-prop}
Let $f_1,\ldots,f_n\in R=k[[x_1,\ldots,x_n]]$ be such that $\AA=R/(f_1,\ldots,f_n)$ is finite-dimensional.
Then the element
$$\bde=\det(\De_j(f_i))\in \AA\ot\AA$$
is equal to the inverse tensor of 
the nondegenerate symmetric bilinear form $(\cdot,\cdot)$ on $\AA$ given by
$$(f,g)=\Res(f\cdot g),$$
where $\Res$
is the Grothendieck residue
\begin{equation}\label{trace-f-def}
\Res(h)=\Res_{k[x]/k}\begin{bmatrix} h(x) dx_1\we\ldots\we dx_n \\ f_1, \ldots,  f_n\end{bmatrix}.
\end{equation}
\end{prop}

\Pf . We give here an argument suggested by one of the referees. In the case $k=\C$ there is also an
analytic proof based on deforming $\w$ to make its critical points nondegenerate (see the proof of
the equality (17) of \cite{KR}).

First of all, let us restate the assertion in a more explicit form. We have to check that
\begin{equation}\label{A-de-eq}
((\cdot,\cdot)\ot\id)(g\ot\bde)=g
\end{equation}
for every $g\in \AA$.
It is known that the bilinear form 
$(\cdot,\cdot)$ is nondegenerate (see \cite[III.9.(R8)]{Hart-RD}), so it is sufficient to check that 
both sides of \eqref{A-de-eq} have the same
pairing with an arbitrary element $h\in \AA$. Thus, we have to check the equality
\begin{equation}\label{traces-eq}
\Res_{k[x,y]/k}\begin{bmatrix} h(x)g(y)\bde(x,y)d\bx\we d\by  \\ 
f_1(x),\ldots,f_n(x),f_1(y),\ldots,f_n(y)\end{bmatrix}=\Res_{k[y]/k}
\begin{bmatrix} h(y)g(y)d\by  \\ f_1(y),\ldots,f_n(y)\end{bmatrix},
\end{equation}
where $d\bx=dx_1\we\ldots\we dx_n$, $d\by=dy_1\we\ldots\we dy_n$.
To do this we will use two standard properties of the Grothendieck residue: Transformation Law and Transitivity 
(see properties (R1) and (R4) in \cite[III.9]{Hart-RD}). 
Denoting $d(\by-\bx)=d(y_1-x_1)\we\ldots\we d(y_n-x_n)$
we can rewrite the left-hand side of \eqref{traces-eq} as
\begin{align*}
&\Res_{k[x,y]/k}\begin{bmatrix} h(x)g(y)\bde(x,y)d\bx\we d\by  \\ 
f_1(x),\ldots,f_n(x),f_1(y),\ldots,f_n(y)\end{bmatrix} \\
&=\Res_{k[x,y]/k}\begin{bmatrix} h(x)g(y)\bde(x,y)d(\by-\bx)\we d\by \\ 
f_1(y)-f_1(x),\ldots,f_n(y)-f_n(x),f_1(y),\ldots,f_n(y)\end{bmatrix} &\text{(Transformation Law)} \\
&=\Res_{k[y]/k}\begin{bmatrix} R(y)g(y)d\by  \\
f_1(y),\ldots,f_n(y)\end{bmatrix}, &\text{(Transitivity)}
\end{align*}
where
$$R(y)=\Res_{k[x,y]/k[y]}\begin{bmatrix} h(x)\bde(x,y)d(\by-\bx)  \\ 
f_1(y)-f_1(x),\ldots,f_n(y)-f_n(x)\end{bmatrix}.$$
Now we observe that $\bde(x,y)$ is the transition
determinant between the systems of parameters given by $(y_1-x_1,\ldots,y_n-x_n)$ and 
$(f_1(y)-f_1(x),\ldots,f_n(y)-f_n(x))$. 
Hence, by the Transformation Law,
we have
$$R(y)=\Res_{k[x,y]/k[y]}\begin{bmatrix} h(x)d(\by-\bx)  \\  y_1-x_1,\ldots,y_n-x_n\end{bmatrix}=h(y).
$$
It remains to substitute this into the above expression for the left-hand side of \eqref{traces-eq}.
\ed

\begin{cor}\label{form-cor} 
The canonical bilinear form \eqref{pair-eq} for $\CC=\MF(R,\w)$ after the identification
\eqref{HH-w-w-eq} 
coincides with 
\begin{equation}\label{form-eq2}
\lan f\ot dx, g\ot dx\ran=(-1)^{{n\choose 2}}\Res(f\cdot g),
\end{equation}
where 
$$\Res(f)=\Res_{k[x]/k}\begin{bmatrix} f(x)\cdot dx_1\we\ldots\we dx_n \\ \pa_1\w,\ldots,\pa_n\w\end{bmatrix}
$$
\end{cor}

\Pf . Since the canonical bilinear form is equal to the inverse of the
Chern character $\ch(\De^{\st})\in H(-\w)\ot H(\w)=H(\w)\ot H(\w)$,
the assertion follows from Propositions \ref{form-prop} and \ref{res-prop} applied to
$f_i=\pa_i\w$.
\ed

Now the general categorical \HRR\ formula \eqref{HRR-cat-eq} specializes to the following
 explicit version for \mf s.

\begin{thm}\label{HRR-thm} (i) For $\bar{E},\bar{F}\in\MF(R,\w)$
we have
\begin{equation}\label{HRR-eq2}
\chi(\Hom_{\w}(\bar{E},\bar{F}))=\dim\Hom^{\zer}_{\w}(\bar{E},\bar{F})-
\dim\Hom^{\one}_{\w}(\bar{E},\bar{F})=\lan \ch(\bar{E}),\ch(\bar{F})\ran,
\end{equation}
where $\lan\cdot,\cdot\ran$ is given by \eqref{form-eq2} and the Chern characters are given by
\eqref{ch-eq2}. 

\noindent
(ii) More generally, for $\a\in\Hom_\w(\bar{E},\bar{E})$ and $\b\in\Hom_\w(\bar{F},\bar{F})$
we have
\begin{equation}\label{cardy-eq}
\sTr_k(m_{\a,\b})=\lan\tau^{\bar{E}}(\a),\tau^{\bar{F}}(\b)\ran,
\end{equation}
where $m_{\a,\b}$ is the endomorphism of $\Hom_\w(\bar{E},\bar{F})$ induced by composing
with $\a$ and $\b$ (see \eqref{m-a-b-eq}).
\end{thm}

\Pf . Part (i) is the particular case of (ii) where $\a$ and $\b$ are the identity morphisms.
To prove (ii) we apply the generalized \HRR\ formula \eqref{cardy-cat-eq}. 
The bilinear form $\lan\cdot,\cdot\ran$ was computed in Corollary \ref{form-cor}. 
Lemma \ref{G-dual-lem} for the trivial group $G$ allows to replace $\tau^{\bar{E}^*}(\a^*)$
with $\tau^{\bar{E}}(\a)$. (To prove \eqref{HRR-eq2} alone it is sufficient 
to use the categorical \HRR\ formula \eqref{HRR-cat-eq} instead of \eqref{cardy-cat-eq}.)
\ed

\begin{rem}\label{odd-rem}
If the number of variables $n$ is odd then $\chi(\Hom_{\w}(\bar{E},\bar{F}))=0$ for any \mf s
$\bar{E}$ and $\bar{F}$. Indeed, this follows from Theorem \ref{HRR-thm},
since in this case $HH_0=0$ so all Chern characters vanish. 
This proves in characteristic zero the conjecture of Hailong Dao 
(see \cite[Conj.\ 3.15]{Dao}). In the $\Z$-graded case a different proof (that works in
arbitrary characteristic) was given by Moore, Piepmeyer, Spiroff and Walker in \cite{MPSW}.
\end{rem}

\begin{rem}
The Chern character of the stabilization $k^{\st}\in\MF(R,\w)$ 
of the residue field $k=R/\m$ (viewed as $R/(\w)$-module),
vanishes for any number of variables $n>0$
(see Proposition \ref{ch-k-st-prop} below).
\end{rem}

\begin{ex}\label{xn-ex} 
Let us illustrate the Theorem for 
$\w=x^n$ (where $n\ge 2$) and the Koszul
\mf s $\bar{E}_i=\{x^i,x^{n-i}\}$, $i=1,\ldots,n-1$. Assume first that $i\ge n/2$. Then the space
$\Hom^1_{\w}(\bar{E}_i,\bar{E}_i)$ is isomorphic to $k[x]/(x^{n-i})\cdot\a_i$, where
$\a_i$ is the odd endomorphism of $\bar{E}_i$ given by
$$\a_i(e_0)=x^{2i-n}e_1,\ \ \a_i(e_1)=-e_0.$$
We have 
$$\tau^{\bar{E}_i}(\a_i)=\sTr(\de'\circ\a)\mod (x^{n-1})=nx^{i-1} \mod(x^{n-1}).$$
For $i\le n/2$ we can get a similar description of the maps $\tau^{\bar{E}_i}$
using the relation $\bar{E}_i\simeq\bar{E}_{n-i}[1]$, where $[1]$ is the change of parity functor.
Notice that $\tau^{\bar{E}_i}$ is injective for every $i$. 
Also, we observe that for $\a\in\Hom^1_\w(\bar{E}_i,\bar{E}_i)$ and
$\b\in\Hom^1_\w(\bar{E}_j,\bar{E}_j)$ 
we have
$$\tau^{\bar{E}_i}(\a)\cdot\tau^{\bar{E}_j}(\b)=0$$ 
in $\AA_\w$ provided $(i,j)\neq (n/2,n/2)$. 
On the other hand, in this case the operator $m_{\a,\b}$ is
nilpotent, so $\sTr_k(m_{\a,\b})=0$ in agreement with
the formula \eqref{cardy-eq}. 
Now let us consider the case $i=j=n/2$ (assuming that $n$ is even). 
Then $\a_i^2=-\id$ and the operator $m_{\a_i,\a_i}$ acts as identity on 
$\Hom^0(\bar{E}_i,\bar{E}_i)\simeq k[x]/(x^i)$ and as $-\id$ on
$\Hom^1(\bar{E}_i,\bar{E}_i)\simeq k[x]/(x^i)$. It follows that
$$\sTr_k(m_{\a_i,\a_i})=2i=n.$$
On the other hand,
$$\lan nx^{i-1},nx^{i-1}\ran=\Res \frac{n^2x^{n-2}dx}{nx^{n-1}}=n,$$
which again agrees with \eqref{cardy-eq}.
\end{ex}

\begin{ex}  Consider the Koszul \mf\ $\bar{E}=\{x;x^2+y^2\}$ of
the $D_4$-singularity $\w=x^3+xy^2$. The Milnor ring for $\w$ is
$$\AA_{\w}=R/(3x^2+y^2, 2xy).$$
The formula \eqref{ab-hom}
shows that $\Hom^0_{\w}(\bar{E},\bar{E})\simeq R/(x,x^2+y^2)=R/(x,y^2)$ is $2$-dimensional,
while $\Hom^1_{\w}(\bar{E},\bar{E})=0$. Hence, 
$$\chi(\bar{E},\bar{E})=2.$$
On the other hand, we have
$$\pa_y\de_E=2y\cdot\iota(e^*) \text{ and } \pa_x\de_E=e\we?+2x\cdot\iota(e^*),$$
where $1,e$ is the standard basis of $\bar{E}$. Hence
$$\ch(\bar{E})=\sTr(\pa_y\de_E\circ\pa_x\de_E)\cdot dx\we dy=
2y\cdot dx\we dy\in\AA_{w}\ot dx\we dy.$$
Hence,
$$\lan\ch(\bar{E}),\ch(\bar{E})\ran=-\Res(4y^2)=\Res\begin{bmatrix} -4y^2\cdot dx\we dy \\ 
(3x^2+y^2), 2xy \end{bmatrix}.
$$
To compute this generalized residue we change the variables to $u=\sqrt{3}x+y$ and 
$v=\sqrt{3}x-y$, and observe that $u^2=(3x^2+y^2)+\sqrt{3}(2xy)$ and 
$v^2=(3x^2+y^2)-\sqrt{3}(2xy)$. Therefore, the above expression is equal to
$$\Res\begin{bmatrix} -4(\frac{u-v}{2})^2\cdot du\we dv \\ 
u^2, v^2\end{bmatrix}=2$$
in agreement with the \HRR\ formula \eqref{HRR-eq2}.
\end{ex}

\subsection{The equivariant case}

Let $G$ be a finite group acting on $R=k[[x_1,\ldots,x_n]]$, and let
$\w\in R$ be a $G$-invariant isolated singularity.

In order to compute the canonical bilinear form on $HH_*(\MF_G(R,\w))$ 
we first calculate the $G\times G$-equivariant Chern character
of the $G\times G$-equivariant stabilized diagonal (see section \ref{sec:G-hoch})
\begin{equation}\label{De-G-st-eq}
\De_{G\times G}^{\st}=\bigoplus_{g\in G}(\id\times g)^*\De_G^{\st}.
\end{equation}
We are going to use the formula \eqref{G-ch-eq} 
for the component $\ch_{G\times G}(\De_{G\times G}^{\st})_{(g_1,g_2)}$ 
with respect to the decomposition
$$
HH_*(\MF_{G\times G}(R^e,\wt{\w}))=
\Bigl(\bigoplus_{(g_1,g_2)\in G\times G}
H(-\w_{g_1})\otimes H(\w_{g_2})
\Bigr)^{G\times G}.
$$
 Note that the element  $(g_1,g_2)$ acts on $\De_{G\times G}^{\st}$
via the isomorphisms
$$
(\id\times h)^*\De_G^{\st}
\rTo{(\id\times h)^*\a_{g_1}} 
(\id\times h)^*(g_1,g_1)^*\De_G^{\st}
\simeq 
(g_1\times g_2)^*
(\id\times g_1hg_2^{-1})^*\De_G^{\st},
$$
where the operator $\a_{g_1}$ is the action of $g_1$ on $\De_G^{\st}$.
Thus, the only non-zero contributions to the 
supertrace may come from the summands corresponding to $h\in G$,
such that $g_1hg_2^{-1}=h$, i.e., $hg_2h^{-1}=g_1$.
Choosing the variables $(x_1,\ldots,x_n)$ in such a way that $g_1$ 
acts by linear transformations preserving $\Span(x_1,\ldots,x_r)$,
and $\Span(x_{r+1},\ldots,x_n)$ is exactly the subspace of $g_1$-invariants
in $\Span(x_1,\ldots,x_n)$, we obtain
\begin{eqnarray}\label{G-str-calc}
&&\ch_{G\times G}(\De_{G\times G}^{\st})_{(g_1,g_2)}=\nonumber
\\
&&\sum_{h\in G:\ hg_2h^{-1}=g_1}
\sTr((\id\times
h)^*[\pa_{y_n}\d\circ\ldots\circ\pa_{y_{r+1}}\d\circ\pa_{x_n}\d\circ\ldots\circ\pa_{x_{r+1}}\d\circ  
g_1]|_{x_1=\ldots=x_r=y_1=\ldots=y_r=0})
\nonumber
\\ 
&&\mod\ (\JJ_{\w_{g_1}}\ot 1,1\ot\JJ_{\w_{g_1}}),
\end{eqnarray}
where the supertrace is computed on the free $R^e$-module $\bigwedge^\bullet(V)\ot R^e$
of the diagonal \mf\ $\De_G^{\st}=(\bigwedge^{\bullet}(V)\ot R^e,\d)$
and $V=\m/\m^2$.
The element $g_1\in G$ acts by the automorphism of the exterior algebra induced by the
action of $g_1$ on $V$.
In particular, $g_1(e_j)=e_j$ for $j>r$, and $g_1$ preserves the subalgebra generated by 
$e_1,\ldots,e_r$.

As in section \ref{sec:calc}
we use the decomposition $\d=\d_{\De}+\d_{\w}$, where
$$\d_{\De}=\iota(\sum_{j=1}^n e_j^*\ot (y_j-x_j)) \text{ and }
\d_{\w}=(\sum_{j=1}^n e_j\ot \wt{\De}_j\w)\we ?$$
and $\wt{\De}_j$ is defined by \eqref{G-phi-eq}.
Let us further split $\d^{\De}$ into two parts $\d_{\De}=\d_{\De}^{\le r}+\d_{\De}^{>r}$, where
$$\d_{\De}^{>r}=\iota(\sum_{j=r+1}^n e_j^*\ot (y_j-x_j)).$$
For $j>r$ we have $\pa_{y_j}\d_{\De}=\pa_{y_j}\d_{\De}^{>r}$ and $\pa_{x_j}\d_{\De}=
\pa_{x_j}\d_{\De}^{>r}$. Hence, the operators $\pa_{y_j}\d$ and $\pa_{x_j}\d$, as well as $g$, 
preserve the filtration 
$$F_p={\bigwedge}^{\ge p}(\bigoplus_{i=1}^r k\cdot e_i)\ot
{\bigwedge}^{\bullet}(\bigoplus_{j=r+1}^n k\cdot e_j)\ot R^e$$
on ${\bigwedge}^{\bullet}(V)\ot R^e$. Hence,
we can pass to the induced endomorphism of the associated graded space, which allows us to
replace $\d=\d_{\De}+\d_{\w}$ with $\d_{\De}^{>r}+\d_{\w}^{>r}$, where
$$\d_{\w}^{>r}=(\sum_{j=r+1}^n e_j\ot \wt{\De}_j\w)\we?.$$
The restriction of $\d_{\De}^{>r}+\d_{\w}^{>r}$ to
$x_1=\ldots=x_r=y_1=\ldots=y_r=0$ coincides with 
the differential for the stabilized diagonal $\De^{\st}_{\w_{g_1}}$ of the potential
$\w_{g_1}=\w|_{x_1=\ldots=x_r=0}$ (tensored with the identity
on $\bigwedge^{\bullet}(\bigoplus_{i=1}^r k\cdot e_i)$). 
Thus, the right-hand side of \eqref{G-str-calc} can be rewritten as
$$\sum_{h\in G:\ hg_2h^{-1}=g_1}
(\id\times h)^*\ch(\De^{\st}_{\w_{g_1}})\cdot\det[\id-g_1;V/V^{g_1}],$$
where the determinant is equal to
the supertrace of $g_1$ acting on $\bigwedge^{\bullet}(\bigoplus_{i=1}^r
k\cdot e_i)$ by the well-known property of the characteristic polynomial. 

This brings us to the following $G$-equivariant version of the formula for the canonical pairing
on the Hochschild homology (cf. Corollary \ref{form-cor})
and of the \HRR\ formula for \mf s\ (cf. Theorem \ref{HRR-thm}).
Recall that we have an isomorphism \eqref{G-hoch-eq}
$$HH_*(\MF_G(R,\w))\simeq \Bigl(\bigoplus_{g\in G} H(\w_g)\Bigr)^G,$$
where $\w_g$ is the restriction of the potential $\w$ to the subspace of $g$-invariants
(we can assume that $G$ acts by linear transformations, see Section
\ref{sec:G-hoch}). 

\begin{thm}\label{G-HRR-thm} 
(i) Let
$$
\lan\cdot,\cdot\ran:HH_*(\MF_G(R,-\w))\ot HH_*(\MF_G(R,\w))\to k
$$
be the canonical bilinear form \eqref{pair-eq} for $\CC=\MF_G(R,\w)$.
Then for
\be
\begin{array}{l}
(h_g)_{g\in G}\in HH_*(\MF_G(R,-\w))=(\bigoplus_{g\in G} H(-\w_g))^G \text{ and } \nonumber\\
(h'_g)_{g\in G}\in HH_*(\MF_G(R,\w))=(\bigoplus_{g\in G} H(\w_g))^G 
\end{array}
\ee
we have
\be \label{G-form-eq}
\lan(h_g),(h'_g)\ran=|G|^{-1}\cdot \sum_{g\in G} c_g\cdot
\lan h_g, h'_g \ran_{\w_g},
\ee
where $\lan\cdot,\cdot\ran_{\w_g}$ is the canonical pairing \eqref{form-eq2} for the potential $\w_g$
and 
$$c_g=\det[\id-g;V/V^g]^{-1},$$
where $V=\m/\m^2$ and $V^g\sub V$ is the subspace of $g$-invariants.

\noindent
(ii) For $\bar{E},\bar{F}\in\MF_G(R,\w)$ 
we have
\be\label{G-HRR-eq}
\chi(\Hom_{\w}(\bar{E},\bar{F})^G)=
|G|^{-1}\cdot\sum_{g\in G}c_g\cdot\lan\ch_G(\bar{E})_{g^{-1}},\ch_G(\bar{F})_g\ran_{\w_g},
\ee
where $\ch_G(\bar{E})_g$ is given by \eqref{G-ch-eq}. 

More generally, for 
$\a\in\Hom_\w(\bar{E},\bar{E})^G$ and $\b\in\Hom_\w(\bar{F},\bar{F})^G$
we have
$$\sTr(m_{\a,\b})=|G|^{-1}\cdot\sum_{g\in G}c_g\cdot
\lan\tau^{\bar{E}}(\a)_{g^{-1}},\tau^{\bar{F}}(\b)_g\ran_{\w_g},$$
where $m_{\a,\b}$ is the endomorphism of $\Hom_\w(\bar{E},\bar{F})^G$ given by \eqref{m-a-b-eq}.
\end{thm}

\Pf . 
(i) Since the canonical bilinear form $\lan\cdot,\cdot\ran$ is equal to the
inverse of the tensor  
$\ch_{G\times G}(\De_{G\times G}^{\st})$ (see \eqref{inv-De-eq}), we have to
check that this tensor 
is equal to the Casimir element of the nondegenerate form \eqref{G-form-eq} on
$\Bigl(\bigoplus_{g\in G} H(\w_g)\Bigr)^G$. 
Now the assertion follows from the similar result in the non-equivariant case 
(see Section \ref{non-equiv-form-sec})
and the calculation preceding the Theorem: we just have to use the fact that
the Casimir element corresponding to the restriction of a $G$-invariant
metric to the subspace of $G$-invariants is given by the 
averaging 
$$\frac{1}{|G|}\sum_{h\in G}(\id\otimes h)T,$$
where $T$ is the Casimir element of the original metric.

 \noindent
(ii) This follows from the generalized \HRR\ formula \eqref{cardy-cat-eq}, 
the equation \eqref{G-form-eq} and Lemma \ref{G-dual-lem}.
\ed

\begin{ex}
Assume that the ground field $k$ contains an $n$th primitive root of unity $\zeta$.
Consider the potential $\w=x^n\in R=k[[x]]$ with the group of symmetries $G=\Z/n$, where an
element $[m]\in \Z/n$ acts by $x\mapsto\zeta^m\cdot x$.
For each $i\in\Z$ let us denote by  $\rho_i$ the character of $G$ given by 
$$\rho_i([m])=\zeta^{mi}.$$
For $i=1,\ldots,n-1$, let us define the $G$-equivariant matrix factorization $\bar{E}_i$ of $x^n$ by
setting $(E_i)^0=\rho_i\ot R$, $(E_i)^1=R$, such that after forgetting the $G$-equivariant structure
we get $\bar{E}_i=\{x^i,x^{n-i}\}$.
Note that in this case $H(\w)^G=0$, while $H(\w_{[m]})^G=H(\w_{[m]})=k$ for $[m]\neq [0]$.
The formula \eqref{G-ch-eq} in this case reduces to
$$\ch_G(\bar{E})_{[m]}=\sTr([m]|_{x=0})$$
for $[m]\neq [0]$.
Thus, we obtain for $[m]\neq [0]$:
$$\ch_G(\rho_a\ot\bar{E}_i)_{[m]}=\rho_{a+i}([m])-\rho_a([m])=\zeta^{am}(\zeta^{mi}-1).$$
On the other hand, $c_{[m]}=(1-\zeta^m)^{-1}$. Thus, we obtain
$$\chi(\bar{E}_i,\rho_a\ot \bar{E}_i)=n^{-1}\cdot
\sum_{m=1}^{n-1}\frac{(\zeta^{-mi}-1)\zeta^{am}(\zeta^{mi}-1)}{1-\zeta^m}.$$
A straightforward calculation allows us to rewrite the right hand side as
$$\sum_{j=0}^{i-1}\de_{[a],[-j]}-\sum_{j=1}^i\de_{[a],[j]}.$$
This agrees with the fact that 
$$\Hom^0_{\w}(\bar{E}_i,\bar{E}_i)=\rho_0\oplus\rho_1\oplus\ldots\oplus\rho_{i-1} \text{ and}$$
$$\Hom^1_{\w}(\bar{E}_i,\bar{E}_i)=\rho_{-1}\oplus\rho_{-2}\oplus\ldots\oplus\rho_{-i}.$$
\end{ex}

\subsection{Boundary-bulk map for the stabilization of the residue field}
\label{k-st-sec}

Here we will compute the Chern character and the boundary-bulk map \eqref{iota-def-eq} for the
stabilization $k^{\st}$ of the residue field $k=R/\mg$.
Recall that if we present $\w$ as
\be\label{w-i-eq}
\w=x_1\w_1+\ldots+x_n\w_n \text{ for }\w_i\in R,
\ee
then $k^{\st}\in\MF(R,\w)$ is the Koszul matrix factorization
$$k^{\st}\{\w_1,\ldots,\w_n;x_1,\ldots,x_n\}=({\bigwedge}^{\bullet}(V)\ot R,\de),$$
where $V=\mg/\mg^2$ and
\be\label{k-st-diff-eq}
\de=(\sum_i e_i\ot\w_i)\we ?+\iota(\sum_i e_i^*\ot x_i)
\ee
with $e_i=x_i\mod \mg^2\in V$.
In the case when $\w$ is preserved by a finite group
of automorphisms $G$, we can equip $k^{\st}$ with a $G$-equivariant structure (see
\eqref{G-k-st}). The key property of $k^{\st}$ is that for any $\bar{E}\in\MF(R,\w)$ 
there is an isomorphism 
\be\label{E-k-eq}
\Hom_{\w}(\bar{E},k^{\st})\simeq (E|_0)^*,
\ee
where $E|_0$ is the restriction of $E$ to the origin (see \cite[Lem.\ 4.2]{Dyck}).
In particular, we obtain an isomorphism of $\Z/2$-graded vector spaces
\be\label{End-k-eq}
H:=\Hom_{\w}(k^{\st},k^{\st})\simeq {\bigwedge}^\bullet\left(\bigoplus_{i=1}^n k\cdot e_i^*\right).
\ee
Let us determine the algebra structure on $H$.

\begin{prop}\label{End-k-prop}  
Assume that $\w\in\mg^2$ and choose
$\w_{ij}\in R$ such that
\be\label{w-i-j-eq}
\w_j=\sum_{i=1}^n x_i\w_{ij} \text{ for } j=1,\ldots,n.
\ee
Then for each $j=1,\ldots,n$, the element
\begin{equation}\label{a-j-eq}
\a_j=-(\sum_i e_i\ot \w_{ij})\we ?+\iota(e_j^*)\in\Homb^1_\w(k^{\st},k^{\st})
\end{equation}
is closed. The cohomology classes $[\a_i]\in H$ generate $H$ as a $k$-algebra
and satisfy the relations
\be\label{Cliff-rel-eq}
[\a_i]\cdot [\a_j]+[\a_j]\cdot[\a_i]=-\w_{ij}(0)-\w_{ji}(0).
\ee
In other words, $H$ is isomorphic to the Clifford algebra associated with the quadratic
form given by the matrix
$(e_i^*,e_j^*)=-\w_{ij}(0)-\w_{ji}(0)$. 
\end{prop}

\Pf . By direct computation we see that $\de\circ\a_j+\a_j\circ\de=0$
and
$$\a_i\circ \a_j+\a_j\circ\a_i=-(\w_{ij}+\w_{ji})\cdot\id.$$
This shows that $\a_j$ is closed.
To deduce \eqref{Cliff-rel-eq} 
we 
combine this with the fact that
$f\cdot\id$ is a coboundary for any $f\in\mg$. Indeed, if $f=x_1f_1+\ldots+x_nf_n$
then 
$$f\cdot\id=[\de^{\st},(\sum_i e_i\ot f_i)\we ?].$$
Thus, the subalgebra $H'\sub H$ generated by the classes $(\a_i)$ is isomorphic to the Clifford algebra.
Since this subalgebra maps bijectively to ${\bigwedge}^\bullet\left(\bigoplus_{i=1}^n k\cdot e_i^*\right)$
under the isomorphism \eqref{End-k-eq}, we conclude that $H'=H$.
\ed

\begin{cor}
If $\w\in\mg^3$ then $H$ is supercommutative and 
\eqref{End-k-eq} is an isomorphism of algebras.
\end{cor}

\begin{rem} In fact, the algebra $H$ is equipped with an $A_{\infty}$-structure.
In the case when $\w\in\mg^3$, the potential $\w$ can be recovered from this $A_{\infty}$-structure
up to a change of variables (see \cite[Thm.\ 7.1]{Ef}).
\end{rem}

Now let us calculate 
the $G$-equivariant Chern character of the stabilization of the residue field
$k^{\st}$ (where $G$ is a finite group of symmetries
of $\w$).

\begin{prop}\label{ch-k-st-prop} 
Let $G$ be a finite group of symmetries of $\w\in R$ and let
 $k^{\st}\in\MF_G(R,\w)$ be the $G$-equivariant
stabilization of $k$ (see \eqref{G-k-st}). Then
$$\ch_G(k^{\st})_g=\begin{cases} \det(\id-g;V), & \text{ if }V^g=0,\\ 0, & \text{{\rm otherwise}}.\end{cases}
$$
\end{prop}

\Pf . We can choose $\w_i$ in \eqref{w-i-eq} so that the differential \eqref{k-st-diff-eq} is
$G$-equivariant (see \ref{sec:G-hoch}).
When $V^g=0$ the formula \eqref{G-ch-eq} gives
$$\ch_G(k^{\st})_g=\sTr(g;{\bigwedge}^{\bullet}(V))=\det(\id-g;V).$$
Now assume that $V^g\neq 0$ and let us choose variables $(x_1,\ldots,x_n)$ in such a
way that the action of $g$ is given by 
$$g(x_1,\ldots,x_n)=(\ell_1,\ldots,\ell_r,x_{r+1},\ldots,x_n),$$
where $\ell_1,\ldots,\ell_r$ are linear forms in $x_1,\ldots,x_r$, and 
$\Span(x_{r+1},\ldots,x_n)$ is exactly the subspace of $g$-invariants in
$\Span(x_1,\ldots,x_n)$. Note that by assumption $r<n$.
By \eqref{G-ch-eq}, we need to prove that in this case
$$\sTr([\pa_n\d\circ\ldots\circ\pa_{r+1}\d\circ g]|_{x_1=\ldots=x_r=0})=0.$$
Consider the filtration $\ldots\supset F_p\supset F_{p+1}\supset\ldots$ on $\bigwedge^{\bullet}(V)\ot R$
with
$$F_p={\bigwedge}^{\ge p}(\bigoplus_{i=1}^r k\cdot e_i)\ot
{\bigwedge}^{\bullet}(\bigoplus_{j=r+1}^n k\cdot e_j)\ot R.$$
By passing to the associated graded space, 
as in the computation of the supertrace in \eqref{G-str-calc},
we can assume that $r=0$. Thus, the problem is reduced to 
the case when $g$ acts trivially, and we have to show
that
$$\sTr_R(\pa_n\d\circ\ldots\circ\pa_1\d)=0$$
for $n>0$. 
Consider
the composition
\be\label{p-decomp}
\pa_n\d\circ\ldots\circ\pa_1\d=(\iota(e_n^*)+p_n\we ?)\circ\ldots\circ(\iota(e_1^*)+p_1\we ?),
\ee
where
\be\label{p-j-eq}
p_j=\sum_{i=1}^n e_i\ot\pa_j\w_i.
\ee
After expanding the right-hand side of \eqref{p-decomp}
only the terms which contain equal amounts of $\iota(e_i^*)$ and  $p_j\we ?$ factors will
contribute to the supertrace. 
Now the assertion follows from Lemma \ref{exterior-tr-lem} below.
\ed

\begin{lem}\label{exterior-tr-lem} 
Let $V$ be a $k$-vector space with the basis $e_1,\ldots,e_n$.
Suppose that we have operators $A_1,\ldots,A_r$ on $K=\bigwedge^\bullet(V)\ot R$,
such that for each $i=1,\ldots,r$, either $A_i=\iota(e_{m}^*)$ for some $m$,
or $A_i=(\sum_{j=1}^n e_j\ot f_{j})\we ?$ for some $f_{1},\ldots,f_{n}\in R$.
Then 
$$\sTr_R(A_1\circ\ldots\circ A_r)=0$$ unless all 
the operators $\iota(e_1^*),\ldots,\iota(e_n^*)$ appear among $A_1,\ldots,A_r$.
\end{lem}

\Pf .
Let $I\sub\{1,\ldots,n\}$ be the set of all $i$ such that $\iota(e_i^*)$ appears among
$A_1,\ldots,A_r$. Consider the decomposition
$$V=V_I\oplus V'_I \text{ where }V_I:=\bigoplus_{i\in I} k\cdot e_i \text{ and }
V'_I=\bigoplus_{j\not\in I} k\cdot e_j.$$
Then each operator $A_i$ preserves the filtration 
$${\bigwedge}^{\ge w}(V'_I)\ot {\bigwedge}^{\bullet}(V_I)\ot R$$
on $\bigwedge^{\bullet}(V)\ot R$.
After passing to the associated graded space, $A_i$ induces an
operator of the form $\id\ot\bar{A}_i$, where $\bar{A}_i$ acts on $\bigwedge^\bullet(V_I)\ot R$.
Thus, we obtain 
$$\sTr_R(A_1\circ\ldots\circ A_r)=\sTr_k(\id;{\bigwedge}^{\bullet}V'_I)\cdot
\sTr_R(\bar{A}_1\circ\ldots\circ\bar{A}_r;{\bigwedge}^{\bullet}V_I\ot R)=0$$
provided $\dim V'_I>0$, i.e., $I$ is a proper subset of $\{1,\ldots,n\}$.
\ed

\begin{rem}
For every $G$-equivariant \mf\ $\bar{E}=(E,\d)$ of $\w$ and a representation $\rho$ of $G$
there is an isomorphism 
\be\label{E-k-rho-eq}
\Hom_{\w}(\bar{E},k^{\st}\ot\rho)^G\simeq \Hom(E|_0,\rho)^G,
\ee
where $E|_0$ is the restriction of $E$ to the origin (see \cite[Lem.\ 4.2]{Dyck}).
The \HRR\ formula \eqref{G-HRR-eq} together with the formula
\eqref{G-ch-eq} and the above Proposition give the following
expression for the Euler characteristic of the left-hand side of \eqref{E-k-rho-eq}:
$$\chi (\Hom_{\w}(\bar{E},k^{\st}\ot\rho)^G)=|G|^{-1}\cdot
\sum_{g, V^g=0}\ch_G(\bar{E})_{g^{-1}}\tr(g;\rho)=
|G|^{-1}\cdot \sum_{g, V^g=0}\sTr(g^{-1};E|_0)\tr(g;\rho).$$
This is compatible with the standard formula for the Euler characteristic of the right-hand side of
\eqref{E-k-rho-eq} because, as we will show,
$$\sTr_k(g;E|_0)=0 \text{ when } V^g\neq 0.$$
Indeed, we can assume that $g$ acts by linear transformations. Furthermore, replacing $G$ by the 
cyclic subgroup
generated by $g$, the \mf\ $\bar{E}$ by its restriction
to the subspace of $g$-invariants and the potential 
$\w$ by $\w_g$, we can assume that $G$ acts trivially on $R$.
Then we have a decomposition 
$$\bar{E}=\bigoplus_{i=1}^N \rho_i\ot \bar{E}_i,$$ 
where $\bar{E}_i$ are 
(non-equivariant) \mf s of $\w_g$ and $\rho_i$ are representations of $G$. Since $\w_g\neq 0$,
the superdimension of 
each $\Z/2$-graded space $\bar{E}_i|_0$ vanishes. Hence, we have
$$\sTr_k(g;E|_0)=\sum_{i=1}^N \tr(g;\rho_i)\cdot\sTr_k(\id;\bar{E}_i|_0)=0.$$
\end{rem}

\

Now let us consider the non-equivariant situation. We are going to calculate the 
boundary-bulk map for $k^{\st}$ assuming that $\w\in\mg^2$.
Recall that by Proposition \ref{End-k-prop},
the algebra $\Hom_\w(k^{\st},k^{\st})$ can be identified with a certain Clifford algebra. 

\begin{prop}\label{Hessian-prop} 
Let $\w\in\mg^2$ and elements $(\w_i)$ and $(\w_{ij})$ are chosen as in 
\eqref{w-i-eq} and \eqref{w-i-j-eq}. Let $[\a_1],\ldots,[\a_n]$ be the generators of the algebra
$\Hom_\w(k^{\st},k^{\st})$ given by \eqref{a-j-eq}. Then 
the boundary-bulk map for $k^{\st}$ is given by
$$\tau^{k^{\st}}([\a_{i_1}]\circ\ldots\circ[\a_{i_r}])=0 \text{ for }r<n$$
and
$$\tau^{k^{\st}}([\a_1]\circ\ldots\circ[\a_n])=
\frac{\Hess(\w)}{\mu}\cdot d\bx \mod\ \JJ_\w\cdot d\bx,$$
where $\Hess(\w)=\det(\pa_i\pa_j\w)$ is the Hessian and 
$\mu=\dim\AA_\w$ is the Milnor number of $\w$.
\end{prop}

\Pf . Recall that
$$\a_j=\iota(e_j^*)-s_j\we?,$$
where 
$$s_j=\sum_{i=1}^n e_i\ot \w_{ij}.$$
Hence, the formula \eqref{tau-eq} gives in our case
\begin{equation}
\begin{array}{l}
(-1)^{{n\choose 2}}\cdot \tau^{k^{\st}}([\a_{i_1}]\circ\ldots\circ[\a_{i_r}])=
\nonumber\\
\sTr_R\left(
(\iota(e_1^*)+p_1\we?)\circ\ldots\circ(\iota(e_n^*)+p_n\we?)\circ
(\iota(e_{i_1}^*)-s_{i_1}\we?)\circ\ldots\circ(\iota(e_{i_r}^*)-s_{i_r}\we?)\right)\cdot d\bx,
\end{array}
\end{equation}
where $p_j$ is given by \eqref{p-j-eq}. By Lemma \ref{exterior-tr-lem}, this expression is zero
for $r<n$. If $(i_1,\ldots,i_r)=(1,\ldots,n)$, we get
\be\label{A-prod-eq}
\begin{array}{l}
(-1)^{{n\choose 2}}\cdot \tau^{k^{\st}}([\a_1]\circ\ldots\circ[\a_n])=\\
\sum_{I\sub\{1,\ldots,n\}}\sTr_R(A_I(p_1,1)\ldots A_I(p_n,n)A_{I^c}(-s_1,1)\ldots 
A_{I^c}(-s_n,n))\cdot
d\bx,
\end{array}
\ee
where $I^c$ denotes the complement of $I$ and
$$A_I(v,i)=\begin{cases} \iota(e_i^*), & i\in I,\\ v\we?, &i\not\in I\end{cases}.$$
Using Lemma \ref{exterior-tr-lem} again we see that we can skew-permute the operators
in the product under the supertrace in \eqref{A-prod-eq}.
Thus, exchanging $A_I(p_i,i)$ with $A_{I^c}(-s_i,i)$ for each $i\in I$ 
produces the factor $(-1)^{|I|}$,
and we get
\be
\begin{array}{l}
\sTr_R(A_I(p_1,1)\ldots A_I(p_n,n)A_{I^c}(-s_1,1)\ldots A_{I^c}(-s_n,n))=\\
\sTr_R\left((v(I,1)\we\ldots\we v(I,n)\we?)\circ\iota(e_1^*)\circ\ldots\circ\iota(e_n^*)\right)=\nonumber\\
(-1)^{{n\choose 2}}\cdot \det(v(I,1),\ldots,v(I,n)),
\end{array}
\ee
where 
$$v(I,i)=\begin{cases}s_i, &i\in I\\ p_i, &i\not\in I\end{cases}.$$
Summing over all subsets $I$ in $\{1,\ldots,n\}$ we obtain
$$\tau^{k^{\st}}([\a_1]\circ\ldots\circ[\a_n])=\det(s_1+p_1,\ldots, s_n+p_n)\cdot d\bx=
\det(\pa_j\w_i+\w_{ij})\cdot d\bx \mod\ \JJ_\w\cdot d\bx.$$
Using \eqref{w-i-eq} and \eqref{w-i-j-eq} we get
$$\pa_j\w=\sum_{i=1}^n x_i(\pa_j\w_i+\w_{ij}).$$
Now the assertion follows from Lemma \ref{res-lem} below applied to $f_i=\pa_i\w$.
\ed

\begin{lem}\label{res-lem} Let $f_1,\ldots,f_n\in \mg\sub R=k[[x_1,\ldots,x_n]]$ be
a regular sequence. Choose $f_{ij}\in R$ such that
\begin{equation}\label{f-i-j-eq}
f_j=\sum_{i=1}^n x_if_{ij}.
\end{equation}
Then 
$$\det(\pa_i f_j)=\mu\cdot\det(f_{ij}) \mod\ (f_1,\ldots,f_n),$$
where $\mu=\dim_k(R/(f_1,\ldots,f_n))$.
\end{lem}

\Pf .
Recall that by the general residue theory (see \cite[III.9.(R8)]{Hart-RD}), 
the invariant pairing 
$$(a, a')=\Res \begin{bmatrix} a\cdot a'\cdot dx_1\we\ldots\we dx_n \\ f_1,\ldots,f_n\end{bmatrix}
$$
on the algebra $\AA=R/(f_1,\ldots,f_n)$ is perfect. 
Furthermore, for any $a\in \AA$, by \cite[III.9.(R6)]{Hart-RD}, one has 
$$(a, \det(\pa_i f_j))=\Res \begin{bmatrix} a\cdot df_1\we\ldots\we df_n \\ f_1,\ldots,f_n
\end{bmatrix}=\Tr_{\AA}(a).$$
For $a$ in the maximal ideal of $\AA$ this expression is zero which implies that
$\det(\pa_i f_j)$ belongs to the socle of $\AA$. Also, setting $a=1$ we get
\begin{equation}\label{Res-J-eq}
\Res \begin{bmatrix} \det(\pa_i f_j)\cdot dx_1\we\ldots\we dx_n \\ f_1,\ldots,f_n\end{bmatrix}=\mu.
\end{equation}
On the other hand, using the adjoint matrix to $(f_{ij})$ one can immediately check that
$\det(f_{ij})x_l$ belongs to $(f_1,\ldots,f_n)$, so $\det(f_{ij})$ belongs to
the socle of $\AA$. 
Also, using \eqref{f-i-j-eq} and the transformation law for the residue (see \cite[III.9.(R1)]{Hart-RD})
we obtain that
\begin{equation}\label{Res-f-i-j-eq}
\Res \begin{bmatrix} \det(f_{ij})\cdot dx_1\we\ldots\we dx_n \\ f_1,\ldots,f_n\end{bmatrix}=1.
\end{equation}
Since the socle is one-dimensional, comparing \eqref{Res-J-eq} with \eqref{Res-f-i-j-eq}
we obtain the required formula.
\ed

\subsection{Graded matrix factorizations}
\label{graded-sec}

Let $L$ be a commutative group with a fixed element $\ell\in L$ such that the quotient
$L/\lan \ell\ran$ is finite.
Assume that the ring $R=k[[x_1,\ldots,x_n]]$ is $L$-graded in such a way that each $x_i$ is homogeneous.
An $L$-graded free $R$-module is a free $R$-module equipped with an $L$-grading such that the
basis elements are homogeneous.

\begin{defn}
For a potential $\w\in R$, homogeneous of degree $2\ell$, an {\it $L$-graded \mf\ } of $\w$
is a pair of finitely generated $L$-graded free $R$-modules $(E^0, E^1)$ equipped with $R$-linear maps
$\de_0:E^0\to E^1$ and $\de_1:E^1\to E^0$, homogeneous of degree $\ell$, 
such that $\de_0\de_1=\w\cdot\id$ and $\de_1\de_0=\w\cdot\id$.
\end{defn}

Equivalently, we can view an $L$-graded \mf\ as a $\Z/2\times L$-graded free $R$-module $E=E^0\oplus E^1$, 
equipped with an endomorphism $\de$ of bidegree $(1,\ell)\in\Z/2\times L$, such that
$\de^2=\w\cdot\id$.

For a pair of $L$-graded \mf s $\bar{E}=(E,\de_E)$ and $\bar{F}=(F,\de_F)$ we define a $\Z$-graded
complex $\Homb_{\w,L}(\bar{E},\bar{F})$ by setting
$$\Homb_{\w,L}(\bar{E},\bar{F})^i=\Hom_{\gr-\Mod_R}^i(E,F(i\cdot\ell)),$$
where $\Hom_{\gr-\Mod_R}^0$ (resp., $\Hom_{\gr-\Mod_R}^1$) is the space of morphisms of 
$\Z/2\times L$-graded $R$-modules of bidegree $(0,0)$ (resp., $(1,0)$).
The differential $d$ on $\Homb_{\w,L}(\bar{E},\bar{F})$ is given by the usual formula \eqref{diff-mf-eq}.
Note that the requirement that $\de_E$ and $\de_F$ have bidegree $(1,\ell)\in\Z/2\times L$ 
implies that $d$ has degree $1$.
In this way we get a dg-category of $L$-graded \mf s.

Consider the finite commutative group $G=L/\lan2\ell\ran$, and let $G^*=\Hom(G,k^*)$
be its dual group (we assume that $k$ contains a primitive root of unity of order $|G|$). 
An $L$-grading on a vector space $V$ induces a natural action of $G^*$ on $V$,
so that $\ga\in G^*$ acts on $V_l$, where $l\in L$, by the scalar multiplication with 
$\ga(l\mod\lan2\ell\ran)\in k^*$.
In particular, we have an action of $G^*$ on $R$ by algebra automorphisms. 
This action preserves $\w$, since $\w$ has degree $2\ell$.

Now suppose we have an $L$-graded \mf\  $\bar{E}=(E,\de_E)$.
Then the $L$-grading on $E$ induces a $G^*$-action, so 
$\bar{E}$ can be viewed as a $G^*$-equivariant \mf\ of $\w$.
Notice that for a pair of $L$-graded \mf s $\bar{E}$ and $\bar{F}$ one
has an equality of $\Z/2$-graded complexes
$$\Homb_{\w,L}(\bar{E},\bar{F})=\Homb_{\w}(\bar{E},\bar{F})^{G^*}.$$

Thus, we can apply Theorem \ref{G-HRR-thm} to calculate 
$$\chi(\bar{E},\bar{F})=\sum_{i\in\Z}(-1)^i\dim H^i(\Homb_{\w,L}(\bar{E},\bar{F}))$$
for $L$-graded \mf s $\bar{E}$ and $\bar{F}$.

{\sc Department of Mathematics, University of Oregon, Eugene, OR 97405}

{\it Email addresses}: apolish@uoregon.edu, vaintrob@uoregon.edu

\end{document}